\documentclass[arxiv,reqno,twoside,a4paper,12pt]{amsart}

\usepackage{tikz}
\usepackage{amsmath, verbatim}
\usepackage{amssymb,amsfonts,mathrsfs,mathtools}
\usepackage[colorlinks=true,linkcolor=blue,citecolor=blue]{hyperref}
\usepackage[driver=pdftex,heightrounded=true,centering]{geometry}
\usepackage{longtable, slashed}
\usepackage{tabularx}
\usepackage{multirow}
\usepackage{caption}
\usepackage{latexsym,enumitem}
\usepackage[all]{xy}

\usepackage[latin2]{inputenc}
\usepackage[mathscr]{eucal}
\usepackage{t1enc, indentfirst}
\usepackage{graphicx, graphics}
\usepackage{pict2e, epic, amssymb }
\usepackage{tikz, pst-node}
\usepackage{tikz-cd, pgfplots}
\usetikzlibrary{arrows}
\usetikzlibrary{calc, patterns}
\usepackage{textcomp}
\pgfplotsset{compat=1.9}
\usepackage{epstopdf}

\usepackage[scaled]{helvet} 
\usepackage{courier} 
\usepackage[mathbf]{euler}

\theoremstyle{plain}
\newtheorem{Th}{Theorem}[section]

\newtheorem{Cor}[Th]{Corollary}

\newtheorem{theorem}[Th]{Theorem}

 \theoremstyle{definition}
\newtheorem{Def}[Th]{Definition}
\newtheorem{defn}[Th]{Definition}

\newtheorem{Rem}[Th]{Remark}
\newtheorem{?}[Th]{Problem}
\newtheorem{Ex}[Th]{Example}
\newtheorem{assum}[Th]{Assumption}

\newcommand{\R}{\mathbb{R}}

\numberwithin{equation}{section}

\definecolor{qqwuqq}{rgb}{0,0,0}

\begin{document}

\title[Heat kernel asymptotics]{Spectral geometry on manifolds with 
fibred boundary metrics II: heat kernel asymptotics}

\author{Mohammad Talebi}
	\address{Universit\"at Oldenburg, Germany}
	\email{mohammad.talebi@uni-oldenburg.de}

	\author{Boris Vertman}
	\address{Universit\"at Oldenburg, Germany}
	\email{boris.vertman@uni-oldenburg.de}
	
	\address{Data availability statement: All data generated or analysed during this study are 
	included in this published article and its supplementary information files.}

\maketitle
\begin{abstract} 
In this paper we continue with the analysis of spectral problems in the setting of complete manifolds with 
fibred boundary metrics, also referred to as $\phi$-metrics, as initiated in our previous work \cite{TGV}. We consider the Hodge Laplacian 
for a $\phi$-metric and construct the corresponding heat kernel as a polyhomogeneous conormal distribution on an appropriate manifold with corners. 
Our discussion is a generalization of an earlier work by Albin and Sher, and provides a fundamental first step towards
analysis of Ray-Singer torsion, eta-invariants and index theorems in the setting.
\end{abstract}

\tableofcontents

\section{Introduction and statement of the main results}

Consider a compact smooth Riemannian manifold $\overline{M}$ with boundary $\partial M$,
which is the total space of a fibration $\phi: \partial M \to B$ over a closed manifold $B$ with the fibre given by 
a closed manifold $F$. Consider a collar neighborhood $[0,\varepsilon]_x \times \partial M$
of the boundary, with a boundary defining function $x$. In the open interior $M$ of such a manifold there 
are various possible complete Riemannian metrics, specified by their structure as $x\to 0$. We shall recall here the main three classes of these 
complete Riemannian metrics.

\subsection{Fibred boundary and scattering metrics}

In this work we are interested in the fibred boundary metrics, 
also called $\phi$-metrics. Ignoring cross-terms for the purpose of a clear exposition, 
these metrics are asymptotically given near the boundary $\partial M$ by 
\begin{equation*}
g_{\phi} = \frac{dx^{2}}{x^{4}}+ \frac{\phi^{*}g_{B}}{x^{2}}+g_{F},
\end{equation*}
where $g_B$ is a Riemannian metric on the base $B$, and $g_F$ is a symmetric bilinear 
form on $\partial M$, restricting to Riemannian metrics on fibres $F$. In case of trivial 
fibres, such a metric is called scattering. A trivial example of a scattering metric is the Euclidean space, 
with the metric written in polar coordinates as $dr^{2} + r^{2}d\theta^{2}$.
After a change of variables $x=r^{-1}$ we obtain
\begin{equation*}
g_{\text{Sc}} = \frac{dx^{2}}{x^{4}} + \frac{d\theta^{2}}{x^{2}}.
\end{equation*}
Such metrics arise naturally in various geometric 
examples. Complete Ricci flat metrics are often $\phi$-metrics. Scattering metrics include 
metrics of locally Euclidean (ALE) manifolds. Products of these spaces with any compact manifold 
provide natural examples of $\phi$-metrics. Furthermore, common classes of gravitational instantons, 
such as the Taub-NUT metrics and reduced 2-monopole moduli space metric, are $\phi$-metrics under 
appropriate coordinate change, cf. \cite[p.2]{HHM}.
\medskip

While there are various approaches to Euclidean scattering theory, a microlocal approach
has been taken by Melrose \cite{Mel-scattering}, where elliptic theory of scattering metrics has been
developed. Elliptic theory of $\phi$-metrics has been studied by Mazzeo and Melrose in \cite{maz-pseudo}.
This work later was generalized to the case of towers of fibrations with so called a-metrics by Grieser and Hunsicker 
\cite{grieser-pseudo}. Elliptic theory of \cite{maz-pseudo} has also been extended by Grieser and Hunsicker \cite{grieser-pseudo2}
to include not fully elliptic operators. Hodge theory of $\phi$-metrics has been developed by Hausel, Hunsicker and Mazzeo 
\cite{HHM}. Index theory (bypassing usual heat operator approach and using adiabatic limit methods instead) in this setting has been addressed by
Leichtnam, Mazzeo and Piazza \cite{LMP}. 

\subsection{Fibred boundary cusp and $b$-metrics}
Fibred (boundary) cusp metrics, also referred to as $d$-metrics 
are conformally equivalent to $\phi$-metrics by a conformal factor $x^2$. 
Ignoring as before the cross-terms for the purpose of a clear exposition, 
these metrics are asymptotically given near the boundary $\partial M$ by 
\begin{equation*}
g_{d} = \frac{dx^{2}}{x^{2}}+x^2 g_{F}+ \phi^*g_{B}.
\end{equation*}
In case of trivial fibres, such a metric is called a $b$-metric. Under the coordinate
change $x= e^{-t}$, a $b$-metric becomes a cylindrical metric $dt^2 + g_B$. Same
change of coordinates turns a genuine fibred cusp metric into $dt^2 + \phi^*g_{B}
+ e^{-2t} g_F$, which is a Q-rank one cusp when $\phi$ is a fibration of tori over a torus. 
Other examples include products of compact manifolds with locally symmetric 
spaces with finite volume hyperbolic cusps. 
\medskip

Elliptic theory of $b$-metrics was pioneered by Melrose \cite{MelATP}. 
Since $\phi$-metrics and fibred cusp metrics differ by a conformal change,
elliptic theory of $\phi$-metrics is suited for fibred cusp metrics as well. Vaillant
\cite{vaillant} has utilized elliptic theory of $\phi$-metrics as well as a microlocal heat kernel construction 
in order to establish an index theorem for fibred cusp metrics. We emphasize that his
heat kernel construction refers to the Hodge Dirac and Hodge Laplacian of a fibred cusp metric, 
not a $\phi$-metric. 

\subsection{Complete edge and conformally compact metrics}
The third class of complete Riemannian metrics on manifolds with fibred boundary, that 
has been of focal relevance in recent geometric analysis developments are complete edge metrics
that by definition are given asymptotically near the boundary $\partial M$ by 
\begin{equation*}
g_e = \frac{dx^{2} + \phi^*g_{B}}{x^{2}} + g_{F}.
\end{equation*}
In case of trivial fibres, such metrics are also called conformally compact with the
classical example being the hyperbolic space $\mathbb{H}^n$. The edge metrics also 
generalize the $b$-metrics that arise as special case of edge metrics with trivial base. 
The significance of edge metrics also lies in their conformal equivalence to the 
incomplete singular wedge metrics. These metrics appeared prominently in the resolution of 
the Calabi-Yau conjecture on Fano manifolds, cf. Donaldson \cite{Donaldson}, Tian \cite{Tia-stability} as well as
Jeffres, Mazzeo, Rubinstein \cite{JMR}.\medskip

Elliptic theory of edge metrics has been developed by Mazzeo \cite{MazEll}. Prior to that, the zero-calculus
containing geometric operators associated to conformally compact metrics, as well as Hodge theory 
have been studied by Mazzeo \cite{Maz-thesis}. Meromorphic extension of the resolvent of conformally
compact spaces is due to Mazzeo and Melrose \cite{MM-mero}. Heat kernel and the (renormalized) Gauss Bonnet index theorem 
on general edge metrics is due to Albin \cite{albin}. Let us also mention the work by Mazzeo and Vertman \cite{maz-ana} on 
analytic torsion and by Vertman \cite{Ver-Ricci} on incomplete wedge spaces, both of which are based on a microlocal
heat kernel construction on wedge manifolds.  

\subsection{Main result and structure of the paper}
While the previous overview of the various classes of metrics on fibred boundary spaces puts our work into 
context, we are dealing here only with fibred boundary $\phi$-metrics and construct the heat kernel 
for finite times microlocally as a polyhomogeneous function on an appropriate heat space blowup.
The construction is similar in spirit to \cite{albin, maz-ana} and \cite{vaillant}.
Our main result is as follows (cf. Corollary \ref{main-corr} for the precise statement)

\begin{theorem}\label{theorem-main}
The heat kernel of the Hodge Laplacian for a $\phi$-metric is a polyhomogeneous 
function on an appropriate heat space blowup for finite time. Same holds for 
Laplace-type operators with normal operator having the same structure as the scalar Laplacian of a
$\phi$-metric, times an identity matrix. This includes for instance square of the spin Dirac operator.
\end{theorem}

The heat kernel construction in the special case of scattering metrics has been outlined in 
the appendix to Sher \cite{Sher-heat}, where the author specifies an initial heat parametrix 
inside the heat calculus of Albin \cite[Theorem 4.3]{albin}. The heat calculus of Albin \cite[Theorem 4.3]{albin}
refers to both scattering and $\phi$-metrics, but concentrates only on the case of complete 
edge metrics in detail. Moreover, \cite{albin} does not deal with the actual heat kernel construction in the $\phi$-case,
as its focus lies on the complete edge case. Our paper closes this gap in the literature. \medskip

The main difference between \cite{Sher-heat, albin} and our presentation here is threefold. First, our initial heat parametrix
construction generalizes the appendix in \cite{Sher-heat} to the more general case of fibered 
boundaries. Second, we do not use the heat calculus of \cite[Theorem 4.3]{albin}, but rather a
simpler version. Third is rather a minor but saddle point, that we work out the Volterra series argument,
which is not addressed in the references above.\medskip

The structure of the paper is as follows. In \S \ref{preliminaries-section} we
introduce the basic geometric preliminaries of $\phi$-metrics, including the 
$\phi$-vector fields and differential $\phi$-operators. Moreover, 
we recall the basic definition of polyhomogeneous functions on manifolds
with corners. In \S \ref{elliptic-section} we review the $\phi$-double space, which arises in the microlocal description of the 
resolvent for the Hodge Laplacian of a $\phi$-metric. That space is built upon in the
heat calculus construction as outlined in the subsequent section \S \ref{heat-section}. 
The heat kernel construction then proceeds in three steps. We first construct the 
heat space and define the heat calculus in \S \ref{space-section}. We then construct
and initial heat parametrix, solving the heat equation to first order in \S \ref{initial-section}.
Our final chapter \S \ref{triple-section} is concerned with the triple space argument, which 
finishes the heat kernel construction and proves the main theorem.

\subsection{Outlook and upcoming work}
In the upcoming work by the first author, the presented heat kernel asymptotics 
together with the low energy resolvent, as constructed in our previous work \cite{TGV} jointly with Daniel Grieser, 
is applied to define the renormalized heat trace and study its asymptotics. This leads to the definition of
a renormalized analytic torsion for $\phi$-manifolds. \medskip

We also intend to apply our analysis to establish an index theorem for
$\phi$-metrics by direct heat kernel methods, following the ideas of \cite{albin} instead of adiabatic arguments of
\cite{LMP}. Note that in view of the heat kernel construction, presented here, Albin \cite[Theorem 6.2]{albin} 
now applies to $\phi$-manifolds and provides an important first step towards an index theorem 
in this setting by heat kernel methods. However, an index theorem will not be a mere
Corollary of our work here: without the image of the Dirac operator being closed, and
without Getzler rescaling techniques, two of the three terms in the formula of \cite[Theorem 6.2]{albin} 
remain mysterious.

\subsubsection*{Acknowledgements} 
The first author acknowledges constructive discussions with
Daniel Grieser and Collin Guillarmou. 
Both authors thank the anonymous referee for careful reading of the manuscript and
valuable suggestions. We also thank university of Oldenburg for financial support and hospitality.

\section{Preliminaries on fibred boundary manifolds and $\phi$-metrics}\label{preliminaries-section}

\subsection{$\phi$-metrics}
Let $\overline{M} = M \cup \partial M$ be a compact smooth manifold with boundary $\partial M$. 
Assume that the boundary $\partial M$ is the total space of the fibration $\phi: \partial M \longrightarrow B$ 
with typical fibre $F$, where both $F$ and base $B$ are smooth compact manifolds without boundary.

\begin{assum}\label{assum-subm}
Throughout the paper we assume that the boundary fibration $(\partial M, g_{\partial M})$ and 
its base $(B,g_B)$ are equipped with Riemannian metrics such that $\phi: (\partial M, g_{\partial M}) 
\to (B, g_B)$ is a Riemannian submersion: we split the tangent bundle $T\partial M$ into the vertical
subbundle $T^V\partial M$ and its orthogonal complement with respect to $g_{\partial M}$ $-$ the horizontal subbundle 
$T^H \partial M$. Then $\phi$ is a Riemannian submersion if the restriction $d\phi: T^H \partial M \to TB$
is an isometry. In this case we may write $g_{\partial M} = g_F + \phi^* g_B$, where $g_F$ equals 
$g_{\partial M}$ on $T^V\partial M$ and vanishes on $T^H\partial M$. \end{assum}

Let $x \in C^{\infty}(\overline{M})$ be a boundary defining function i.e $x \geq 0$, $\partial M = \{x = 0\}$ and $dx \neq 0$ at $\partial M$. 
By the collar neighborhood theorem there always exists a collar $[0,\varepsilon)_x \times \partial M \subset \overline{M}$
of the boundary. Replacing $g_F + \phi^* g_B$ by $\frac{\phi^*g_B}{x^2} +  g_F$, for any constant positive value of $x$, still defines a Riemannian submersion.
We can now define $\phi$-metrics.

\begin{defn}\label{def-phi-metric}
We call a Riemannian metric $g_{\phi}$ on the open interior $M$ a $\phi$-metric, 
if in the collar neighborhood $ \mathscr{U} = ( 0,\varepsilon)_x \times \partial M$ of the boundary, 
$g_{\phi}$ takes the following form
\begin{equation}
g_{\phi} \restriction \mathscr{U} = \frac{dx^{2}}{x^4} + \frac{\phi^{*}g_{B}}{x^{2}} + g_{F} + h =: g_{0} + h,
\end{equation}
where $h$ is a higher order term, i.e. satisfies 
$\vert h \vert_{g_{0}} = O(x)$ as $x \to 0$.
\end{defn}

The Assumption \ref{assum-subm} is used in order to deduce the structure \eqref{hodge laplace}
for the Hodge Laplacian, which in turn is essential in the second step of the heat kernel construction
in \eqref{NHfd}.

\begin{Ex}\label{example-euclidean} The Euclidean plane $\R^m$ is a particular example of a $\phi$-manifold
(a scattering manifold) with the fibre $F$ being a single point and the base $B= \mathbb{S}^{m-1}$.
Indeed, choosing polar coordinates $(r,\theta)$ and writing $x = 1/r$, the Euclidean metric on $\R^m$
can be written away from the origin as 
$$
g_{\R^m} = \frac{dx^2}{x^4} + \frac{g_{\mathbb{S}^{m-1}}}{x^2}.
$$
\end{Ex}

\subsection{$\phi$-vector fields and differential $\phi$-operators}
 
 \begin{defn}\label{Phi-vector-def}
 The $\phi$-vector fields $\mathcal{V}_{\phi}\equiv \mathcal{V}_{\phi}(M)$
 are by definition smooth vector fields over $\overline{M}$, tangent to the fibres of $\phi$ and such that for 
 $x \in C^{\infty}(\overline{M})$ we have $Vx \in x^{2}C^{\infty}(\overline{M})$. Any $\phi$-vector field
 is then locally generated near $\partial M$ by
 $$
 x^{2}\frac{\partial}{\partial x}, x\frac{\partial}{\partial y_{i}}, \frac{\partial}{\partial z_{j}},
 $$
where $\{x,y_i,z_j\}$ are local coordinates on $[0,\varepsilon) \times \partial M$ with $y = \{y_i\}_i$
being local coordinates on the base $B$, lifted to $\partial M$ and extended to 
 $[0,\varepsilon) \times \partial M$, and $z=\{z_j\}_j$ restricting to local coordinates on the fibres $F$.
We introduce the so called $\phi$-tangent space by requiring $\mathcal{V}_{\phi}(M)$
to be its smooth sections
\begin{equation*}
C^\infty(\overline{M}, {}^\phi TM) = \mathcal{V}_{\phi}(M)
= C^\infty(\overline{M})\text{-span}\left\langle x^{2}\frac{\partial}{\partial x}, 
x\frac{\partial}{\partial y_{i}}, \frac{\partial}{\partial z_{j}}\right\rangle.
\end{equation*}
The dual ${}^\phi T^*M$ , the so-called $\phi$-cotangent space, satisfies
\begin{equation*}
C^\infty(\overline{M}, {}^\phi T^*M) = C^\infty(\overline{M})\text{-span}\left\langle \frac{dx}{x^{2}}, \frac{dy_{i}}{x},dz_{j}\right\rangle.
\end{equation*}
 \end{defn}

The space of $\phi$-vector fields $\mathcal{V}_{\phi}(M)$ has a Lie-Algebra structure 
and is a $C^{\infty}(\overline{M})$-module. Therefore one may introduce $\textup{Diff}_{\phi}^{k}(M)$ as an graded 
algebra. Explicitly, $ P \in \textup{Diff}_{\phi}^{k}(M)$ if it is a $k$-th order differential operator 
in the open interior $M$ of the following structure near the boundary
\begin{equation}\label{p-phi}
 P = \sum_{\vert \alpha \vert + \vert \beta \vert + q \leq k} 
 P_{\alpha,\beta,q}(x,y,z)(x^{2}D_{x})^{q}(xD_{y})^\beta D_{z}^{\alpha},
 \end{equation}
with coefficients $P_{\alpha,\beta,q}\in C^\infty(\overline{M})$ smooth up to the boundary. 
Its $\phi$-symbol $\sigma_{\phi,k}(P)$ is then locally given for any cotangent vector 
$(\xi, \eta, \zeta) \in {}^\phi T^*M$ over the base point $(x,y,z) \in \overline{M}$ by 
the homogeneous polynomial
\begin{equation*}
 \sigma_{\phi,k}(P) (x,y,z; \xi, \eta, \zeta)= \sum_{\vert \alpha \vert + \vert \beta \vert + q = k} 
 P_{\alpha,\beta,q}(x,y,z)\xi^{q}\eta^\beta \zeta^{\alpha}.
 \end{equation*}
We say that $P$ is $\phi$-elliptic if $\sigma_{\phi,k}(P)$ is
invertible off the zero-section of ${}^\phi T^*M$.
Writing $P^{k}({}^\phi T^*M)$ for the space of 
homogeneus polynomial of degree $k$ on the fibres of ${}^\phi T^*M$, the 
$\phi$-symbol map defines a short exact sequence
\begin{equation}
0 \longrightarrow \textup{Diff}_{\phi}^{k-1}(M) \hookrightarrow \textup{Diff}_{\phi}^{k}(M) 
\overset{\sigma_{\phi,k} \ }{\longrightarrow} P^{k}({}^\phi T^*M) \longrightarrow 0.
\end{equation}
Same constructions extend to case of differential operators acting on sections of 
a flat vector bundle $(L,\nabla)$ over $\overline{M}$, compactly supported in the interior $M$. 
In that case, coefficients $P_{\alpha,\beta,q}$ are smooth sections of the endomorphism bundle 
$\textup{End}(L)$ and each derivative $X \in \mathcal{V}_{\phi}$ in \eqref{p-phi} 
is replaced by $\nabla_X$. We write $\textup{Diff}_{\phi}^{*}(M,L)$ for differential $\phi$-operators acting 
on sections, and the notion of ellipticity carries over verbatim.

\subsection{Hodge-Laplacian of a $\phi$ metric}

The Hodge Laplacian $\Delta_\phi$ is an element of $\textup{Diff}_{\phi}^{2}(M,\Lambda^{*} {}^\phi T^*M)$
and this section is devoted to writing out its explicit structure, following Hausel, Hunsicker and Mazzeo \cite[\S 5.3.2]{HHM}.
\medskip

Assume for the moment that the fibration $\phi$ is trivial, so that 
we can identify $\partial M \cong B \times F$. The exact model $\phi$-metric 
is given in this case  in a collar neighborhood $\mathscr{U}= (0,\varepsilon) \times \partial M$ of the boundary by
$$
g_{0} \restriction \mathscr{U}  = \frac{dx^{2}}{x^{4}} + \frac{g_B}{x^2} + g_F,
$$ 
where in this model case $g_F$ is assumed to be constant along $B$.
We write $b= \dim B$. In this case the Laplace Beltrami operator of $(M,g_0)$ is of the following explicit form in a 
collar neighborhood $\mathscr{U}$ of the boundary
\begin{equation}\label{D-Phi}
\Delta_{\phi} \restriction \mathscr{U} = -x^{4}\partial_{x}^{2} + x^{2}\Delta_{B} + \Delta_{F}  - (2-b)x^{3}\partial_{x}.
\end{equation}
where  $\Delta_B$ is the Laplace Beltrami operator of $(B,g_{B})$ and 
$\Delta_F$ is the family of Laplace Beltrami operator on the fibres $(F,g_{F})$. We want to explain in what way the 
structure of the Hodge Laplacian acting on differential forms for a general $\phi$-metric admits a similar structure as above.
\medskip

Under the Assumption \ref{assum-subm} we may split $T \partial M$ and its dual 
$T^* \partial M$ orthogonally with respect to $\phi^{*}g_{B} + g_{F}$ 
into vertical and horizontal parts. Writing $\mathcal{V}$ for the canonical vertical bundle, and 
$\phi^*TB$ for the horizontal bundle as in Definition \ref{def-phi-metric} (we also write $\phi^*T^*B$
for its dual), we obtain
\begin{equation}\label{splitting1}
T \partial M = \phi^*TB \oplus \mathcal{V}, \quad T^* \partial M = \phi^*T^*B \oplus \mathcal{V}^*.
\end{equation}
This splitting induces an orthogonal splitting of the
$\phi-$cotangent bundle ${}^\phi T^*M$ in the collar neighborhood $\mathscr{U}$ of the boundary
\begin{equation}\label{splitting2}
{}^\phi T^*M \restriction \mathscr{U} = \textup{span} \left\{ \frac{dx}{x^2}\right\} \oplus  x^{-1} \phi^*T^*B \oplus \mathcal{V}^*.
\end{equation}
Let us assume that the higher order term $h\equiv 0$ for the moment. 
With respect to the corresponding decomposition of $\Lambda^{*} {}^\phi T^*M$ over $\mathscr{U}$ we compute
for the exterior derivative, cf. \cite[\S 5.3.2]{HHM} 
\begin{equation}\label{hodge de rham}
D_\phi = x^{2}D_{x} + x \mathbb{A} + D_{F} + xD_{B} 
-x^{2}\mathscr{R}.
\end{equation}
We shall now explain the individual terms in \eqref{hodge de rham}.
The term $x^{2}D_{x}$ acts for any section $\omega$ of the bundle $\Lambda^\ell \left( x^{-1} \phi^*T^*B\right) \oplus \Lambda \mathcal{V}^*$ as follows
$$
\left( x^{2}D_{x}\right)  \omega = \frac{dx}{x^2} \wedge \left( x^{2}\partial_{x} \right) \omega, \quad 
x^{2}D_{x} \left( \frac{dx}{x^2} \wedge \omega \right) = -\left( x^{2}\partial_{x} \right) \omega.
$$
The term $\mathbb{A}$ is given by $A + A^{*}$, where $A$ is a $0-$th order differential operator, acting for any
section $\omega$ of $\Lambda^\ell \left( x^{-1} \phi^*T^*B\right) \oplus \Lambda \mathcal{V}^*$ by
$$
A \omega = -\ell \cdot  \frac{dx}{x^2} \wedge \omega, \quad 
A \left( \frac{dx}{x^2} \wedge \omega \right) = (b-\ell) \cdot \omega.
$$
The term $D_{F} = d_{F} + d_{F}^{*}$, acts as the Gauss Bonnet operator on the
$\Lambda \mathcal{V}^*$ component, cf. the first displayed equation in \cite[p. 527]{HHM}. The term 
$D_B$ is given by
\begin{equation}\label{DB}
D_{B} =  \left( d_{B} - \mathbb{I} \right) +  \left( d_{B} - \mathbb{I} \right)^*. 
\end{equation} 
Here, $\mathbb{I}$ is the second fundamental form, $d_B$ is the lift of the exterior derivative on $B$ to $\partial M$
plus the action of the derivative in the $B$-direction on 
the $\mathcal{V}^*-$components of the form, cf. the second displayed equation in \cite[p. 527]{HHM}.
Finally, $\mathscr{R} = R + R^{*}$, where $R$ is the curvature of 
the Riemannian submersion $\phi$. We can now take the square of $D_\phi$ to compute the Hodge Laplacian
over $\mathscr{U}$
\begin{equation}\label{hodge laplace}\begin{split}
\Delta_{\phi} = D^2_\phi = - (x^{2}D_{x})^2 + D^2_{F} + x^2D^2_{B}
+Q, \ Q \in x\cdot \textup{Diff}^{^2}_{\phi}(\mathscr{U}, \Lambda^{*} {}^\phi T^*\mathscr{U}).
\end{split}\end{equation}
Note that $D^2_B$ equals $\Delta_B$ (the Hodge Laplacian on $B$) up to additional higher order terms in
$x\cdot \textup{Diff}^{^2}_{\phi}(\mathscr{U}, \Lambda^{*} {}^\phi T^*\mathscr{U})$. 
For the general higher order term $h$ with $|h|_{g_0}=O(x)$ as $x\to 0$, the arguments carry over up to higher order and the 
statement \eqref{hodge laplace} is still true. This replaces \eqref{D-Phi} in the general case.

\begin{Rem}\label{dirac-squares} The arguments of this paper apply to a more general class of 
Laplace-type operators in $\textup{Diff}_{\phi}^{2}(M,L)$, provided the structure 
of \eqref{D-Phi} (times the identity matrix) holds in the collar neighborhood $\mathscr{U}$ up to
higher order terms. This includes squares of geometric Dirac operators, such as the 
spin Dirac operator. \end{Rem}

\subsection{Polyhomogeneous functions on manifolds with corners}\label{phg-section}

The contribution of this paper is a heat kernel construction for the 
Hodge Laplacian $\Delta_{\phi}$ such that the heat kernel lifts to a polyhomogeneous distribution 
on an appropriate manifold with corners. In this subsection we provide brief definitions of manifolds with corners,
polyhomogeneous functions on and maps between these spaces. We refer the reader to
\cite{MelATP} and \cite{basics} for more details.
\begin{Def}
An $n$-dimensional compact manifold $X$ with corners is by definition locally 
modelled near each $p \in X$ diffeomorphically by $(\mathbb{R}^+)^{k}\times\mathbb{R}^{N-k}$ for some $k \in \mathbb{N}_0$, where 
we write $\R^+=[0,\infty)$. The index $k$ is called the codimension of $p$. A boundary face of $X$
is the closure of a connected component of the set of points of codimension $1$. We assume that all boundary faces are
embedded, i.e. each boundary hypersurface $H \subset X$ is given by $\{\rho_{H} = 0\}$ 
for some boundary defining function $\rho_{H} \in C^{\infty}(X)$, where 
$d_{\rho_{H}}\neq 0$ and $\rho_{H} \geq 0$. A corner is  the closure of a connected component of the set of points of codimension at least two. 
\end{Def}

\noindent We now define polyhomogeneous functions on 
manifolds with corners. Polyhomogeneous sections valued in vector bundles over 
manifolds with corners are defined analogously. 

\begin{defn}\label{phg}
Let $X$ be a manifold with corners and $\{(H_i,\rho_i)\}_{i=1}^N$ an enumeration 
of its (embedded) boundaries with the corresponding defining functions. For any multi-index $\beta = (b_1,
\ldots, b_N)\in \mathbb{C}^N$ we write $\rho^{\beta} = \rho_1^{b_1} \ldots \rho_N^{b_N}$.  Denote by 
$\mathcal{V}_b(X)$ the space of smooth vector fields on $X$ which lie
tangent to all boundary faces. An index set 
$E_i = \{(\gamma,p)\} \subset {\mathbb C} \times {\mathbb{N}_0}$ 
satisfies the following hypotheses:

\begin{enumerate}
\item the real parts $\textup{Re}(\gamma)$ accumulate only at $+\infty$,
\item for each $\gamma$ there exists $P_{\gamma}\in \mathbb{N}_0$, such 
that $p \leq P_\gamma$ for all $(\gamma, p) \in E_{i}$.
\item if $(\gamma,p) \in E_i$, then $(\gamma+j,p') \in E_i$ for all $j \in {\mathbb{N}_0}$ and $0 \leq p' \leq p$. 
\end{enumerate}
An index family $E = (E_1, \ldots, E_N)$ is an $N$-tuple of index sets. 
Finally, we say that a smooth function $\omega$ on the interior of $X$ is polyhomogeneous  
with index family $E$, we write $\omega\in \mathscr{A}_{\textup{phg}}^E(X)$, 
if near each $H_i$, 
\[
\omega \sim \sum_{(\gamma,p) \in E_i} a_{\gamma,p} \rho_i^{\gamma} (\log \rho_i)^p, \ 
\textup{as} \ \rho_i\to 0,
\]
with coefficients $a_{\gamma,p}$ polyhomogeneous on $H_i$ with index $E_j$
at any intersection $H_i\cap H_j$ of hypersurfaces. We require the asymptotic expansion to be
preserved under repeated application of $\mathcal{V}_b(X)$. Since $\dim H_i < \dim X$,
this is an inductive definition in the dimension of $X$.
\end{defn}

\noindent There are following classes of morphisms between manifolds with corners.

\begin{Def}\label{b-maps-def} 
Let $X$ and $X'$ be two manifolds with corners and $f:X \longrightarrow X'$
a smooth map, i.e. in each local chart $f$ can be extended to a smooth map between open 
(in $\mathbb{R}^n$) domains containing the charts. 
\begin{itemize}
\item[(1)] $f$ is called a b-map if for any collection $\{\rho_j\}_{j\in J}, \{\rho'_i\}_{i\in I}$ 
of boundary defining functions on $X$ and $X'$, respectively, there exist
non-negative numbers $\{\alpha_{ij}\}_{ij}$ such that for any index $i\in I$
\begin{equation}\label{aij}
f^{*}\rho_{i}' = a_{i}\prod_{j \in J}\rho_{j}^{\alpha_{ij}},\quad 0 < a_{i} \in C^{\infty}(X).
\end{equation}
We define b-tangent bundles as follows: in any local chart $(\R^+)^{k}_{(x)}\times\mathbb{R}^{n-k}_{(y)}$
with local coordinates $(x)=\{x_1,\cdots, x_k\}$ and $(y) = \{y_1, \cdots, y_{n-k}\}$ the 
b-tangent bundle is defined by the spanning sections being tangent to the boundaries of the chart, i.e.
$$
x_{1}\partial_{x_{1}},\cdots,x_{k}\partial_{x_{k}},
\partial_{y_{1}}, \cdots, \partial_{y_{n-k}}.
$$
If $f$ is a b-map, restriction of the total differential $df$ to vector fields that
are tangent to the boundaries of the chart above, is valued in vector fields that
are tangent to the boundaries of the chart again. Hence, the total differential of
a b-map induces 
$${}^b df: {}^bTX \to {}^bTX'.$$

\item[(2)] A b-map $f$ is called b-submersion if its total differential $df$
induces a surjective map ${}^b df$ between b-tangent bundles. \medskip

\item[(3)] A b-submersion is called b-fibration if for each $j\in J$ there is at most one $i \in I$ such that $\alpha_{ij} \neq 0$.
The numbers $\alpha_{ij}$ are defined by \eqref{aij}.
\end{itemize}
\end{Def}

The significance of these maps is that polyhomogeneous functions pull back to polyhomogeneous
functions under b-maps. Moreover, push forward under b-fibrations of densities with coefficients given by polyhomogeneous functions 
is again a density with polyhomogeneous coefficients, see \cite{Mel-push} for the explicit statement of the pullback and the pushforward theorems
and the explicit form of the index families.

\section{Construction of the $\phi$-double space $M^2_\phi$}\label{elliptic-section}

The material of this section is drawn from \cite{maz-pseudo}, where a calculus of pseudo-differential
operators on $M$, containing parametrices of elliptic $P \in \textup{Diff}_{\phi}^{m}(M,L)$,
is developed. The Schwartz kernels of these parametrices are built of polyhomogeneous conormal 
distributions on a certain manifold with corners, with the usual singularity along the diagonal. In this 
section we only provide the definition of that manifold with corners rather than explaining other elements
of the $\phi$-calculus. \medskip

Since $\overline{M}$ is a manifold with boundary, $\overline{M}^2$ has two boundary hypersurfaces, intersecting at $B = \partial M \times \partial M 
\subset \overline{M}^2$. Consider near $B \subset \overline{M}^2$ local coordinates $(x,y,z) , (x', y', z')$, which are just two copies of the local 
coordinates introduced in Definition \ref{Phi-vector-def}. In these coordinates, $B= \{x=x'=0\}$. 

\subsection{b-double space $M^2_b := [\overline{M}^2, B]$}
\begin{figure}[h]    
\centering
\begin{tikzpicture}[scale = 1]
    \draw[->] (0,0) -- (0,2);
    \draw[->] (0,0) -- (2,0);
    \node at (-0.2,2.2) {$x$};
    \node at (2.2,-0.2) {$x'$};

    \draw[<-] (2.7,1) -- node[above] {$\beta_b$} (4.3,1);

   \begin{scope}[shift={(5,0)}]
        \draw[-](1,0) --(2,0);
        \draw[-](0,1)--(0,2);
         \draw  (1,0) arc (0:90:1);
          \node at (-0.3,1.7) {$lf$};
           \node at (2.2,-0.2) {$rf$};
           \node at (1,1) {$ff$};
 \end{scope}
 \end{tikzpicture}
\caption{The b-double space $M^2_b$.}
\label{b-double space}
\end{figure}
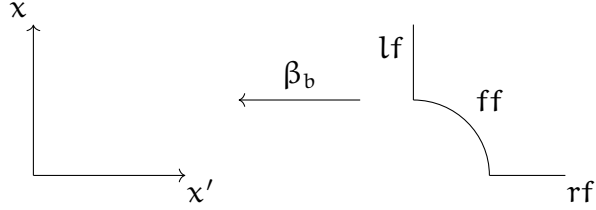

The blowup $[\overline{M}^2, B]$ is defined as the disjoint union of $M^2 \backslash B$ 
with the interior spherical normal bundle of $B \subset \overline{M}^2$, under an identification explained in 
\cite{MelATP}, cf. also \cite{basics}. The blowup $[\overline{M}^2, B]$ is equipped with the minimal differential 
structure such that smooth functions in the interior of $\overline{M}^2$ and polar coordinates around $B$ are smooth. 
The interior spherical normal bundle of $B$ defines a new boundary hypersurface of $[\overline{M}^2, B]$, the front face ff,
in addition to the previous boundary faces $\{x = 0\}$ (the right face) and $\{x'= 0\}$ (the left face).
The front face is itself a quarter circle fibration over $B$. $M^2_b = [\overline{M}^2, B]$ is equipped with the obvious blowdown 
map $\beta_{b}: M_{b}^{2} \to \overline{M}^2$ and is illustrated in Figure \ref{b-double space}.

We now introduce projective coordinates in a neighborhood of the front face ff,
that are very convenient in the computations, but are not globally defined over all of ff.
Near the corner intersection of ff and rf, away from lf, we have projective coordinates
\begin{equation}\label{b-proj-coord1}
s=\frac{x}{x'}, \, y, \, z, \, x', \, y', \, z',
\end{equation}
where $s$ is a defining function of rf and $x'$ a defining function of ff.
Near the corner intersection of ff and lf, away from rf, we get projective coordinates
by interchanging the roles of $x$ and $x'$. The pullback by the blowdown map is locally simply a 
change of coordinates. \medskip

\subsection{$\phi$-double space $M^2_\phi$}
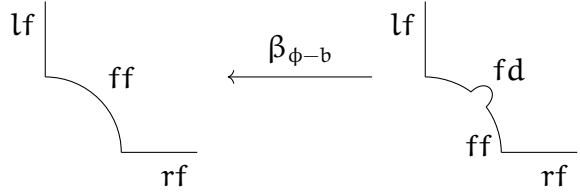
\begin{figure}[h]
\centering

\begin{tikzpicture}[scale = 1]
     \draw[-](1,0) --(2,0);
     \draw[-](0,1)--(0,2);
      \draw  (1,0) arc (0:90:1);
       \node at (-0.3,1.7) {$lf$};
         \node at (1.7,-0.3) {$rf$};
         \node at (1,1) {$ff$};

  \draw[<-] (2.4,1)-- node[above] {$\beta_{\phi-b}$} (4.3,1);

   \begin{scope}[shift={(5,0)}]
        \draw[-](1,0) --(2,0);
        \draw[-](0,1)--(0,2);
         \draw (0,1).. controls (0.3,1) and (0.55,0.85) .. (0.6,0.8);
             \draw (1,0).. controls (1,0.3) and (0.85,0.55) .. (0.8,0.6);
            \draw (0.6,0.8).. controls (0.8,1) and (1,0.8)..(0.8,0.6);
          \node at (1.1,1.1) {$fd$};
         \node at (-0.3,1.7) {$lf$};
         \node at (1.7,-0.3) {$rf$};
         \node at (0.7,0.1) {$ff$};

    \end{scope}

\end{tikzpicture}
\caption{The $\phi$-double space $M_{\phi}^{2}$}
\label{Phi-double space}
\end{figure}

We now blow up the fibre diagonal $\phi \subset M^2_b$ 
\begin{equation*}
\phi :=\{(h,h') \in B = \partial M \times \partial M ; \phi(h) = \phi(h')\},
\end{equation*}
lifted up the b-double space, where $\phi: \partial M \to B$ is the fibration of the boundary. Its lift $\beta_b^{*}(\phi) \subset M^2_b$
is given in the local coordinates \eqref{b-proj-coord1} by $\{y=y',s=1, x'=0\}$. We define the $\phi$-double space by
a similar procedure as above, by blowing up $\beta_b^{*}(\phi)$
\begin{equation*}
M_{\phi}^{2} := [ M_{b}^{2} ; \beta_b^{*}(\phi)], \quad \beta_{\phi-b}: M_{\phi}^{2} \longrightarrow M_{b}^{2}.
\end{equation*}
The new boundary hypersurface is denoted by fd and the blowup is illustrated 
in Figure \ref{Phi-double space}. We also define the full blowdown map
$\beta_{\phi} := \beta_{b} \circ \beta_{\phi-b}: M_{\phi}^{2} \to \overline{M}^2$. 

In local coordinates we have the following projective coordinates near fd
(here the roles of $x$ and $x'$ can be interchanged, leading to equivalent 
projective coordinate systems)
\begin{equation}\label{b-proj-coord2}
S:= \frac{s-1}{x'}, \, U:= \frac{y-y'}{x'}, \, x', \, y', \, z, \, z',
\end{equation}
where $x'$ is the defining function of fd, and the rest of the original 
front face ff lies in the limit $|(S,U)| \to \infty$. As before, the pullback by the 
blowdown map is locally simply a change of coordinates. 

\begin{Rem}
The centrality of the $\phi$-double space $M^2_{\phi}$ stems from the 
fact that Schwartz kernels of parametrices to elliptic differential $\phi-$operators,
studied in the $\phi$-calculus by Mazzeo and Melrose \cite{maz-pseudo}, lift to 
polyhomogeneous conormal distributions on $M^2_{\phi}$ with 
conormal singularity along
the lifted diagonal. 
\end{Rem}

\section{Outline of the heat kernel construction for a $\phi$-metric}\label{heat-section}

Let $\Delta_\phi$ be the unique self-adjoint extension of the Hodge Laplacian on 
the $\phi-$manifold $M$ with fibred boundary $\partial M$ and a $\phi$-metric $g_\phi$. 
As noted in Remark \ref{dirac-squares}, we can replace $\Delta_\phi$ by any 
Laplace-type operator in $\textup{Diff}_{\phi}^{2}(M,L)$, provided the structure 
of \eqref{D-Phi} (times the identity matrix) holds in the collar neighborhood $\mathscr{U}$ up to
higher order terms. This includes squares of geometric Dirac operators, such as the 
spin Dirac operator. \medskip

The heat operator of $\Delta_\phi$ is denoted by $e^{-t\Delta_\phi}$ and solves 
for any given section $\omega_0$ of $\Lambda^* {}^\phi T^*M$ in the domain of $\Delta_\phi$ by definition the homogeneous heat problem 
\begin{equation} \begin{split}
(\partial_t + \Delta_\phi) \omega (t,p)= 0, & \quad (t,p) \in [0,\infty) \times M, \\
\omega(0,p) = \omega_0(p), & \quad p \in M,
\end{split} \end{equation}
with $\omega =e^{-t\Delta_\phi}\omega_0$, a section of $\Lambda^* {}^\phi T^*M$ for any fixed $t \in [0,\infty)$. 
The heat operator is an integral operator
\begin{equation} \label{eqn:hk-on-functions}
e^{-t\Delta_\phi}\omega_0(p) = \int_M H\left( t, p, \widetilde{p} \right) \omega_0(\widetilde{p}) 
\textup{dvol}_{g_\phi} (\widetilde{p}),
\end{equation}
with the heat kernel $H(t, \cdot, \cdot)$ being a section of $ \Lambda^* {}^\phi T^*M \boxtimes \bigl(\Lambda^* {}^\phi T^*M \bigr)^*$
for any $t \in [0,\infty)$. Here, $\boxtimes$ denotes the external tensor product of vector bundles. \medskip

Consider local coordinates $(t, (x,y,z), (\widetilde{x}, \widetilde{y}, \widetilde{z}))$ near the 
highest codimension corner in the heat space $M^2_h$, where $(x,y,z)$ and $(\widetilde{x}, \widetilde{y}, \widetilde{z})$ 
are two copies of the local coordinates on $M$ near the boundary $\partial M$, as before. 
Then the heat kernel $H$ has a non-uniform behaviour at the diagonal $D=\{(0,p,p) \mid p \in \overline{M}\}$ and at the submanifold
\begin{align*}
A = \{ (t, (0,y,z), (0, \widetilde{y}, \widetilde{z}))\in M^2_h : y= \widetilde{y} \}.
\end{align*}
The asymptotic behaviour of $H$ near the submanifolds $D$ and $A$ of $M^2_h$ 
is conveniently studied using the blowup procedure of \S \ref{elliptic-section}. We 
proceed in the remainder of this paper with the following 3 steps: \medskip

\noindent \underline{\emph{Step 1:}} We construct the "heat blowup space" $HM_{\phi}$ by an additional 
blowup in $[0,\infty) \times M_{\phi}^{2}$, where $M_{\phi}^{2}$ is the $\phi$-double space
introduced in \S \ref{elliptic-section}. More specifically, we lift the diagonal $D \in M^2_h$ to 
$[0,\infty) \times M_{\phi}^{2}$ and blow it up, treating $\sqrt{t}$ as a smooth variable\footnote{Note that this is not the parabolic blowup! Namely, in the parabolic blowup,  $\sqrt{t}$ is not smooth away from the blowup. 
However, for the purpose of heat kernel construction, both 
procedures lead to the same result. This is because 
we consider Schwartz kernels that vanish to infinite order as $t \to 0$ away from the blowup.} (i.e. we extend $t \in (0,1)$ to $[0,\infty]$ and consider the smooth structure of functions smooth in $\sqrt{t}$). We then may define the heat calculus of smoothing operators 
with Schwartz kernels that lift to polyhomogeneous functions on 
the heat blowup space $HM_{\phi}$. \medskip

\noindent \underline{\emph{Step 2:}} We obtain an initial
 parametrix for $H$ inside the
heat calculus, solving the heat equation up to first order. This requires us to 
lift the heat equation to $HM_{\phi}$ and to solve the resulting equations 
(normal problems) at the various boundary faces of the heat blowup space. 
\medskip

\noindent \underline{\emph{Step 3:}} The exact heat kernel is then obtained by a Volterra series argument, which requires
the triple space construction and the composition formula of the final section \S \ref{triple-section}.

\section{Step 1: Construction of the heat blowup space}\label{space-section}

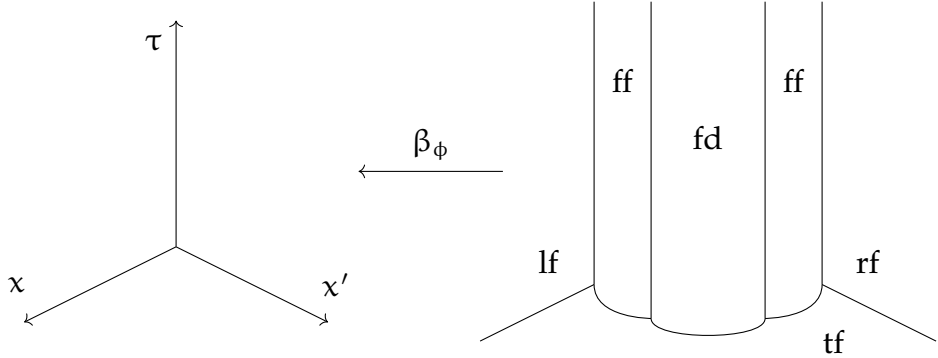
\begin{figure}[h]
\centering
\begin{tikzpicture}[scale = 1]
     \draw[->](0,0) --(2,-1);
     \draw[->](0,0)--(-2,-1);
     \draw[->](0,0)--(0,3);
\node at (2.1,-0.5) {$x'$};
           \node at (-2.1,-0.5) {$x$};
           \node at (-0.3,2.7) {$ \tau$};
 \draw[<-] (2.4,1) -- node[above] {$\beta_{\phi}$} (4.3,1);
\begin{scope}[shift={(7,1)}, scale = 1.5]
        \draw[-](1,-1) --(2,-1.5);
     \draw[-](-1,-1)--(-2,-1.5);
      \draw[-](1,-1) --(1,1.5);
      \draw[-](-1,-1)--(-1,1.5);
      \draw[-](-0.5,-1.3)--(-0.5,1.5);
      \draw[-](0.5,-1.3)--(0.5,1.5);

            \node at (0,0.3) {$\text{fd}$};
         \node at (0.75,0.8) {$\text{ff}$};
         \node at (-0.75,0.8) {$\text{ff}$};
              \node at (1.4,-0.8) {$\text{rf}$};
              \node at (1.1,-1.5) {$\text{tf}$};
         \node at (-1.4,-0.8) {$\text{lf}$};
         \draw (1,-1).. controls (1,-1.3) and (0.5,-1.3)..(0.5,-1.3);
           \draw (-0.5,-1.3).. controls (-0.5,-1.5) and (0.5,-1.5)..(0.5,-1.3);
            \draw (-1,-1).. controls (-1,-1.3) and (-0.5,-1.3)..(-0.5,-1.3);        
\end{scope}
\end{tikzpicture}
\caption{Intermediate heat blowup space $[0,\infty) \times M_{\phi}^{2}$.}
\label{intermediate-blowup}
\end{figure}

Consider the $\phi$-double space $M_{\phi}^{2}$ with the blowdown map 
$\beta_{\phi}: M_{\phi}^{2} \to \overline{M}^2$. We obtain an intermediate heat blowup space
by taking its product with the time axis $[0,\infty)$. We treat the square root of the 
time variable $\tau := \sqrt{t}$ as a smooth variable. The resulting intermediate 
heat blowup space is illustrated in Figure \ref{intermediate-blowup}. \medskip

Let us explain the abbreviations for the boundary hypersurfaces in the 
intermediate heat blowup space: ff stands for front face, fd $-$ the fibre diagonal, 
lf and rf $-$ the left and right faces, respectively, and finally tf stands for the 
temporal face. In view of \eqref{b-proj-coord1} and \eqref{b-proj-coord2}, the projective coordinates on 
$[0,\infty) \times M_{\phi}^{2}$ near the various boundary hypersurfaces are as follows. 

\subsubsection*{Projective coordinates near the intersection of rf with ff, away from fd }

In view of the projective coordinates \eqref{b-proj-coord1} on the $\phi$ double space, 
we have the following projective coordinates, which are valid uniformly up to tf, away from 
an open neighborhood of fd.

\begin{equation}\label{proj-coord1}
s=\frac{x}{x'}, \, y, \, z, \, x', \, y', \, z', \tau = \sqrt{t}.
\end{equation}
In these coordinates, $s$ is a defining function of rf, $\tau$ is a defining 
function of tf, $x'$ is a defining function of ff. Interchanging the roles of $x$ and $x'$
yields projective coordinates near the intersection of lf with ff, where $s'= x' / x$ is a 
defining function of lf. 

\subsubsection*{Projective coordinates near the fd up to tf}

In view of the projective coordinates \eqref{b-proj-coord2} on the $\phi$ double space, 
we have the following projective coordinates, which are valid in open neighborhood of fd
uniformly up to tf, away from lf and rf. Here the roles of $x$ and $x'$ can be interchanged, leading to an 
equivalent system of coordinates.

\begin{equation}\label{proj-coord2}
S= \frac{s-1}{x'}, \, U = \frac{y-y'}{x'}, \, z, \, x', \, y', \, z' , \, \tau = \sqrt{t}.
\end{equation}
In these coordinates, $\tau$ is a defining 
function of tf, $x'$ is a defining function of fd,
and ff lies in the limit $|(S,U)| \to \infty$. 

\subsubsection*{Heat blowup space as a blowup of temporal diagonal} 

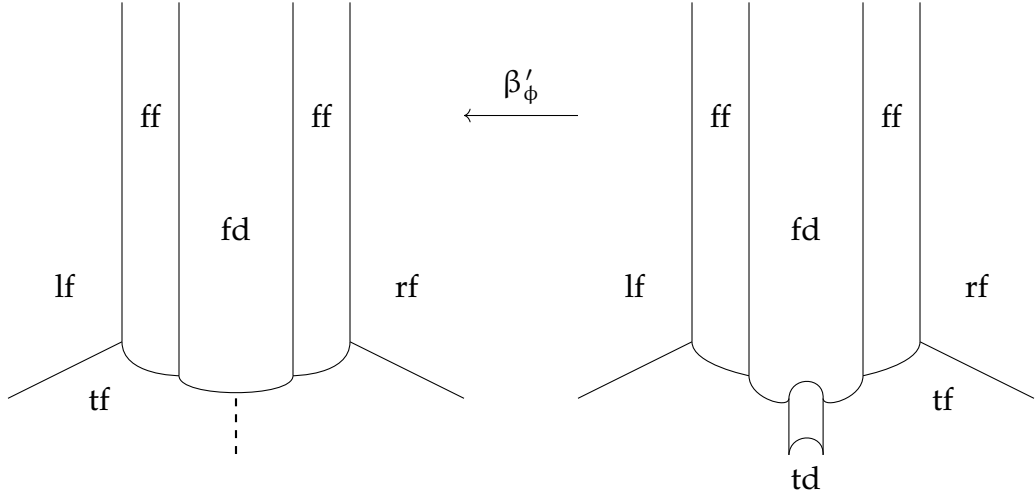
\begin{figure}[h]
\centering

\begin{tikzpicture}[scale = 1.5]
  \draw[-](1,-1) --(2,-1.5);
     \draw[-](-1,-1)--(-2,-1.5);
      \draw[-](1,-1) --(1,2);
      \draw[-](-1,-1)--(-1,2);
      \draw[-](-0.5,-1.3)--(-0.5,2);
      \draw[-](0.5,-1.3)--(0.5,2);
      
            \node at (0,0) {$\text{fd}$};
         \node at (0.75,1) {$\text{ff}$};
         \node at (-0.75,1) {$\text{ff}$};
           \node at (-1.5,-0.5) {$\text{lf}$};
          \node at (1.5,-0.5) {$\text{rf}$};
           \node at (-1.2,-1.5) {$\text{tf}$};
             
         \draw (1,-1).. controls (1,-1.3) and (0.5,-1.3)..(0.5,-1.3);
           \draw (-0.5,-1.3).. controls (-0.5,-1.5) and (0.5,-1.5)..(0.5,-1.3);
            \draw (-1,-1).. controls (-1,-1.3) and (-0.5,-1.3)..(-0.5,-1.3);        
         
      \draw[thick, dashed](0,-1.5)--(0,-2);
          
    \draw[<-] (2,1) -- node[above] {$\beta'_{\phi}$} (3,1);

   \begin{scope}[shift={(5,0)}]
     \draw[-](1,-1) --(2,-1.5);
     \draw[-](-1,-1)--(-2,-1.5);
      \draw[-](1,-1) --(1,2);
      \draw[-](-1,-1)--(-1,2);
      \draw[-](-0.5,-1.3)--(-0.5,2);
      \draw[-](0.5,-1.3)--(0.5,2);
      
          \node at (0,0) {$\text{fd}$};
         \node at (0.75,1) {$\text{ff}$};
         \node at (-0.75,1) {$\text{ff}$};
           \node at (-1.5,-0.5) {$\text{lf}$};
          \node at (1.5,-0.5) {$\text{rf}$};
           \node at (1.2,-1.5) {$\text{tf}$};
           
         \draw (1,-1).. controls (1,-1.2) and (0.5,-1.3)..(0.5,-1.3);
  
            \draw (-1,-1).. controls (-1,-1.2) and (-0.5,-1.3)..(-0.5,-1.3);        
           
       \draw[-] (0.15,-1.5)--(0.15,-2);
        \draw[-](-0.15,-1.5)--(-0.15,-2);
 \draw (-0.15,-1.5).. controls (-0.14,-1.3) and (0.16,-1.3)..(0.15,-1.5);
 \draw (-0.15,-2).. controls (-0.14,-1.8) and (0.16,-1.8)..(0.15,-2);

           \draw (0.15,-1.5)..controls (0.15,-1.6) and (0.5,-1.5)..(0.5,-1.3);
            \draw (-0.15,-1.5)..controls (-0.15,-1.6) and (-0.5,-1.5)..(-0.5,-1.3);

          \node at (0,-2.2) {$\text{td}$};
\end{scope}

\end{tikzpicture}

\caption{The heat blowup space $HM_{\phi}$.}
\label{heat-space}
\end{figure}

The final heat blowup space $HM_{\phi}$ is obtained by blowing up the
lift of the diagonal $D=\{(0,p,p) \mid p \in \overline{M}\}$ to $[0,\infty) \times M_{\phi}^{2}$. 
In the local coordinates in an open neighborhood $\mathscr{U}_{\textup{fd}}$ of fd, its lift is given by 
\begin{align}
\beta_\phi^{*}(D) \restriction \mathscr{U}_{\textup{fd}} = 
\{S=0,U=0,z=z',\tau = 0\}.
\end{align}
The heat blowup space $HM_{\phi}$ is then defined as a parabolic blowup 
\begin{equation}
HM_{\phi} := \left[\, [0,\infty) \times M_{\phi}^{2}, \beta_\phi^{*}(D) \, \right].
\end{equation}
This blowup space is illustrated in \ref{heat-space}. \ \medskip

The full blow down map is defined by 
$$
\beta := \beta_{\phi} \circ \beta'_{\phi}: HM_{\phi}
\to [0,\infty) \times \overline{M}^2.
$$
Let us now describe the resulting heat blowup space in projective coordinates.
The previous coordinate systems \eqref{proj-coord1} and \eqref{proj-coord2}
remain valid away from an open neighborhood of the new boundary face td. 
Near td we have up to the intersection with fd the following projective coordinates
(with respect to the notation of \eqref{proj-coord2})
\begin{equation}\label{proj-coord3}
\mathcal{S} = \frac{S}{\tau} = \frac{x-x'}{(x')^2 \sqrt{t}}, \, 
\mathscr{U} = \frac{U}{\tau} = \frac{y-y'}{x'\sqrt{t}}, \, 
\mathcal{Z} = \frac{z-z'}{\sqrt{t}}, \, 
x',\, y', \, z', \,\tau = \sqrt{t}.
\end{equation}
In these coordinates, $x'$ is a defining function of fd,
 $\tau$ is a
defining function of td, and tf lies in the limit $|(\mathcal{S}, \mathscr{U}, \mathcal{Z})| \to \infty$.
In these coordinates, the roles of $x$ and $x'$ may be
 interchanged freely.
The pullback by the blowdown map $\beta$ is locally simply
 a change between standard and projective 
coordinates \eqref{proj-coord1}, \eqref{proj-coord2} and \eqref{proj-coord2}.
\medskip

We conclude the heat blow up space construction by singling out a class of 
polyhomogeneous functions on it, that define the "heat calculus" in 
our setting.

\begin{Def}[Heat calculus] We write $m= \dim M$. 
Consider a vector bundle $L$ over $\overline{M}$ and the 
projection $\pi: [0,\infty) \times \overline{M}^2 \to \overline{M}^2$ onto the second component. Define the pullback bundle over
$HM_{\phi}$ by
\begin{equation}\label{E}
E:= (\pi \circ \beta)^*  \Bigl( L\, \boxtimes \, L^* \Bigr).
\end{equation}
For any $a,\ell \in \mathbb{R}$, the space $\mathscr{H}^{a,\ell}_{\phi}(M,L)$ is defined as the space of linear operators 
$A$ acting on smooth sections of $L$ that are compactly supported in $M$
$$
A: C^\infty_0 (M, L) \to C^\infty_0 (M, L),
$$
with Schwartz kernels $K_{A}$ lifting to polyhomogeneus sections
 $\beta^{*}K_{A}$ of $E$ such that
for any defining functions $\rho_{\textup{fd}}$ and $\rho_{\textup{td}}$ of fd and td, respectively, we have
\begin{equation*}
\beta^{*}K_{A} = \rho_{fd}^{-3+a}\cdot \rho_{td}^{-m+\ell}G_A.
\end{equation*}
Here $G_A$ is polyhomogeneous, smooth at fd and td,
and vanishing to infinite order at lf, rf, ff and tf. We further define 
$$
\mathscr{H}_\phi^{\infty}(M,L) := \bigcap_{a,\ell, \in \mathbb{R}} 
\mathscr{H}_\phi^{a,\ell}(M,L).
$$
Our main example is $L = \Lambda^{*} {}^\phi T^{*} M$. We shall abbreviate $\Lambda^*_\phi := \Lambda^{*} {}^\phi T^{*} M$
and denote the corresponding heat calculus by $\mathscr{H}_\phi^{a,\ell}(M,\Lambda^*_\phi)$.
\end{Def}

\begin{Ex}
Consider the Example \ref{example-euclidean} and the scalar Euclidean heat kernel $H$ on $\R^m$. 
We want to explain that $H$ lies inside the heat calculus defined above, more precisely 
$H \in \mathscr{H}^{3,0}_{\phi}(M,\Lambda^0_\phi)$. Write 
for any pair $p,p' \in \R^m$
$$
p =: \| p \| \eta =: \eta / x, \quad p' =: \| p' \| \eta' =:  \eta' / x',
$$
where $\eta, \eta' \in \mathbb{S}^{m-1}$. For simplicity, let us assume $\eta = \eta'$. We obtain
from the standard formula for the scalar Euclidean heat kernel, using the projective coordinates \eqref{proj-coord3},
where $\tau$ and $x'$ are defining functions of td and fd, respectively
$$
H(t,p,p') = \frac{1}{(4 \pi)^{\frac{m}{2}}}\, \tau^{-m}\, \exp \left( - \frac{|\mathcal{S}|}{4 |(1+x' \mathcal{S})|} \,  \|  \eta \|^2 \right).
$$
This corresponds to the claim $H \in \mathscr{H}^{3,0}_{\phi}(M,\Lambda^0_\phi)$. Indeed, 
$H$ is smooth at $\{x'=0\}$, i.e. at fd. $H$ behaves as $\tau^{-m}$ at td. It vanishes to infinite order 
at ff, i.e. as $|S| \to \infty$. The general case $\eta \neq \eta'$ is similar, albeit with more cumbersome computations.
\end{Ex}

\begin{Rem} 
We point out that our Schwartz kernels are not multiplied with half-densities here,
which is common in many other references. This simplifies our presentation here, but 
leads to some shifts in the asymptotics later on, when we study compositions of the
Schwartz kernels in  \S \ref{triple-section}. 
\end{Rem}

\section{Step 2: Construction of an initial heat kernel parametrix}\label{initial-section}

We construct an initial heat parametrix for $\Delta_\phi$ by solving the heat equation, lifted to 
$HM_{\phi}$, to leading order at fd and td. The solutions of the heat equation 
at fd and td can be extended off these boundary faces with any power
of the respective defining functions. The correct powers are determined by studying the
lift of the delta distribution at fd and td as well.

\subsubsection*{Solving the heat equation near fd}\label{lifts-subsection}

Let us consider the relevant geometric quantities written in projective coordinates near fd.
Recall the projective coordinates \eqref{proj-coord2}, that are valid near fd, away from td.
We compute in these coordinates
\begin{equation}\label{lifts-derivatives-fd}
\beta^{*}\partial_x = (x')^{-2} \partial_{S}, \quad \beta^{*}\partial_y = (x')^{-1} \partial_{U},  \quad \beta^{*}\partial_z = \partial_z.
\end{equation}
Let us point out that fd is the total space
of fibration over $B$ with fibres $\R \times \R^b \times F^2$. Here, $y' \in B$ denotes the 
base point of the fibration, $(S,U,z,z') \in \R \times \R^b \times F^2$ coordinates on the fibres.
Recall the formula \eqref{hodge laplace} for the Hodge Laplacian $\Delta_\phi$ near the boundary. 
Note that the term $D^2_B$ in \eqref{hodge laplace} equals $\Delta_B$ (the Hodge Laplacian on $B$) plus terms in 
$x\cdot \textup{Diff}^{^2}_{\phi}(\mathscr{U})$. Thus, in view of \eqref{lifts-derivatives-fd} we compute
\begin{align}\label{normal-laplace-fd}
\beta^{*}\Delta_{\phi} \restriction \textup{fd} =
\Delta_{S, U, y'} + \Delta_{F,y'},
\end{align}
where $\Delta_{S,U,y}$ is the Euclidean Hodge Laplacian on $\R^{b+1} \cong \R \times T_{y'} B$ with Euclidean coordinates $(S,U)$,
where the scalar product on $T_{y'} B$ is defined with respect to $g_B(y')$. The second summand, $\Delta_{F,y'}$ is given by $D^2_{F_{y'}}$ 
in the notation of \eqref{hodge laplace}, which is a Laplace type operator acting on sections of $\Lambda \mathcal{V}^*$. We conclude
\begin{align*}
\beta^{*}\left(t\left(\partial_t + \Delta_{\phi}\right)\right) \restriction \textup{fd} = 
\frac{1}{2} \tau \partial_\tau + \tau^2 \left(\Delta_{S, U, y'} + \Delta_{F,y'} \right) =: \mathcal{L}_{\textup{fd}}.
\end{align*}
Note that the parameter $y'$ simply indicates the base point of the fibration fd, and 
for each fixed base point, the equation $\mathcal{L}_{\textup{fd}} u = 0$ is an partial 
differential equation on the fibres of fd. A solution to $\mathcal{L}_{\textup{fd}} u = 0$ is given by
\begin{equation}\label{NHfd}
N_{\textup{fd}}(H)(\tau, S,U,z,z'; y'):= 
H_{\Delta_{S,U,y'}}(\tau, S, U,0) H_{ \Delta_{F,y'}}(\tau, z,z'),
\end{equation}
where $H_{\Delta_{S,U,y'}}$ and $H_{ \Delta_{F,y'}}$ are the heat kernels of 
$\Delta_{S,U,y'}$ and $\Delta_{F,y'}$, respectively. We extend it off the front face fd as 
\begin{equation}\label{H_0}
\beta^* H'_0 := N_{\textup{fd}}(H)(\tau, S,U,z,z'; y') \, \psi(x),
\end{equation}
with a cutoff function $\psi \in C^\infty_0[0,\infty)$ with compact support in $[0,2\varepsilon)$ for some $\varepsilon \in (0,1)$ such that 
$\psi \equiv 1$ on $[0,\varepsilon]$. Then $H'_0$ solves the heat equation to first order at fd, i.e.
$$
\beta^{*}\left(t\left(\partial_t + \Delta_{\phi}\right) H'_0 \right) \restriction \textup{fd} = 0.
$$

\subsubsection*{Solving the heat equation near $td \cup fd$} 

Our heat parametrix $H'_0$ does not solve the 
heat equation to any order at td. Here we explain the standard procedure
how $H'_0$ is corrected to provide a heat parametrix, solving the heat equation to higher order at td as well. 
Recall, near td we have up to the intersection with fd the coordinates \eqref{proj-coord3}
\begin{equation}
\mathcal{S} = \frac{S}{\tau} = \frac{x-x'}{(x')^2 \sqrt{t}}, \, 
\mathscr{U} = \frac{U}{\tau} = \frac{y-y'}{x'\sqrt{t}}, \, 
\mathcal{Z} = \frac{z-z'}{\sqrt{t}}, \, 
x',\, y', \, z', \,\tau = \sqrt{t}.
\end{equation}
In these coordinates, $x'$ is a defining function of fd,
 $\tau$ is a
defining function of td, and tf lies in the limit $|(\mathcal{S}, \mathscr{U}, \mathcal{Z})| \to \infty$.
In these coordinates, the roles of $x$ and $x'$ may be interchanged freely.
We fix coordinates near $td \cap fd$. In these coordinates the individual partial derivatives are written as follows
\begin{equation}\label{lifts-derivatives-td}
\begin{split}
\beta^{*}\partial_x = \frac{1}{\tau x'^2} \partial_{\mathcal{S}}, \quad
\beta^{*}\partial_y = \frac{1}{\tau x'} \partial_{\mathscr{U}}, \quad
\beta^{*}\partial_z = \frac{1}{\tau} \partial_{\mathcal{Z}}.
\end{split}
\end{equation}
Let us point out that td is the total space
of fibration over $\overline{M}$ with fibres $\R \times \R^b \times \R^f$. Here, $p'=(x',y',z') \in \overline{M}$ denotes the 
base point of the fibration, $(\mathcal{S}, \mathscr{U}, \mathcal{Z}) \in \R \times \R^b \times \R^f$ coordinates on the fibres.
In view of \eqref{lifts-derivatives-td} and \eqref{hodge laplace} we compute
\begin{align*}
\beta^{*} t\Delta_{\phi} \restriction \textup{td} = \Delta_{\, \mathcal{S}, \mathscr{U}, \mathcal{Z}, x',y',z'},
\end{align*}
where on the right hand side we have for each fixed $(x',y',z')$ the Hodge Laplacian on 
$\R \times \R^b \times \R^f \cong \R \times T_{y'}B \times \mathcal{V}_{(y'\!, \, z')}$, 
with Euclidean coordinates $(\mathcal{S}, \mathscr{U}, \mathcal{Z})$, defined with respect to the metric $g_B(y')$ at the $T_{y'}B \cong \R^b$
component, and the metric $g_F(y',z')$ at the $\mathcal{V}_{(y'\!, \, z')} \cong \R^f$ component. Away from fd, this is
simply the Hodge Laplacian on $T_{p'}M \cong \R^m$ with respect to the metric $g_\phi$ on $T_{p'}M$. Denote the corresponding 
heat kernel, evaluated at $\tau = 1$, by $N_{\textup{td}}(H)$. We extend it off td, using the cutoff function $\psi$ as above, and define
\begin{equation}\label{H_00}
\beta^* H''_0 := \tau^{-m} N_{\textup{td}}(H)(\tau, \mathcal{S}, \, \mathscr{U}, \, \mathcal{Z}, \, x'\!,\, y'\!,\, z') \, \psi(\tau).
\end{equation}
By construction, $H''_0$ solves the heat equation to first order at td
$$
\beta^{*}\left(t\left(\partial_t + \Delta_{\phi}\right) H''_0 \right) \restriction \textup{td} = O(\tau^{-m+1}), 
\quad \textup{as} \ \tau \to 0.
$$
Since by construction $\beta^*H'_0$ and $\beta^* H''_0$ agree on $td \cap fd$, there 
exists a well-defined operator $H_0 \in \mathscr{H}_\phi^{\, 3,0}(M, \Lambda^*_\phi)$, solving the 
heat equation to first order at fd and td. 

\subsubsection*{Constructing an initial heat parametrix} 

In the final step below, one improves the parametrix $H_0$
to solve away the error at td to infinite order. This is done by the usual argument, which is outlined in various cases, cf. Melrose \cite{MelATP},
and Grieser \cite{notes} as a basic reference, as well as Albin \cite{albin}, Mazzeo and Vertman 
\cite{maz-ana} for the same argument in different geometric settings. This defines a new
heat parametrix, still denoted by $H_0$ in $\mathscr{H}_\phi^{\, 3,0}(M, \Lambda^*_\phi)$ such that 
$$
(\partial_t + \Delta_\phi) H_0 \in \mathscr{H}_\phi^{\, 4,\infty}(M, \Lambda^*_\phi).
$$
Same arguments apply when $\Delta_\phi$ is replaced the square of a geometric 
Dirac operator, such as the spin Dirac operator, or more generally by any 
Laplace-type operator in $\textup{Diff}_{\phi}^{2}(M,L)$, provided the structure 
of \eqref{D-Phi} (times the identity matrix) holds in the collar neighborhood $\mathscr{U}$ up to
higher order terms. This yields the same structure as in \eqref{normal-laplace-fd}, 
so that the rest of the argument applies. This proves the following result.

\begin{Th}\label{initial-parametrix-th}
There exists an initial heat parametrix $H_{0} \in \mathscr{H}_\phi^{\, 3,0}(M, \Lambda^*_\phi)$ for the Hodge 
Laplacian $\Delta_\phi$, solving the 
heat equation to first order at fd, and to infinite order at td, i.e. 
$(\partial_t + \Delta_\phi) H_0 \in \mathscr{H}_\phi^{\, 4,\infty}(M, \Lambda^*_\phi)$.
The restriction of $H_0$ to fd is given by $N_{\textup{fd}}(H)$. The 
leading term in the asymptotic expansion of $H_0$ at td is given by $N_{\textup{td}}(H)$. \medskip

If $\Delta_\phi$ is replaced by a Laplace-type operator in $\textup{Diff}_{\phi}^{2}(M,L)$ 
such that \eqref{D-Phi} (times the identity matrix) holds in the collar neighborhood $\mathscr{U}$ up to
higher order terms, the same statement holds with $ \Lambda^*_\phi \equiv  \Lambda^{*} {}^\phi T^{*} M$
replaced by $L$. 
\end{Th}

\section{Step 3: Triple space construction and composition of operators}\label{triple-section}

In this section we use the initial heat parametrix $H_0$ in Theorem \ref{initial-parametrix-th}
to construct the exact heat kernel as a polyhomogeneous section of $E$ in \eqref{E} on the heat space $HM_\phi$.
The construction is based on the following composition result, which is the main technical 
result of this section.

\begin{Th}(\textbf{Composition Theorem})\label{compth}
 Assume that, $A \in \mathscr{H}_\phi^{a,\ell}(M, \Lambda^*_\phi)$ and 
$B \in \mathscr{H}_\phi^{a',\infty}(M, \Lambda^*_\phi)$. We denote the corresponding Schwartz kernels of $A$ and
 $B$ by $K_A$ and $K_B$, respectively. Then the composition $A \circ B$ with the Schwartz kernel
 given by 
\begin{align}\label{prop comp} 
  K_{A\circ B}(t,p,p') := \int_{M}\int_0^t K_A(t-t',p,p'') \cdot K_B(t',p'',p') dt' \, \textup{dvol}_{g_\phi}(p''),
\end{align}
is well defined and $A \circ B \in \mathscr{H}_\phi^{a+a',\infty}(M, \Lambda^*_\phi)$.
Same statement holds for $\Lambda^*_\phi$ replaced by any vector bundle $L$ over $\overline{M}$.
\end{Th}

\noindent The dot in \eqref{prop comp} is defined fibrewise as follows
\begin{equation}\label{dot-def} \begin{split}
\cdot : \Bigl(  \Lambda^*_{\phi, p}  \boxtimes (\Lambda^*_{\phi,p''})^* \Bigr) 
\times \Bigl( \Lambda^*_{\phi,p''} \boxtimes  (\Lambda^*_{\phi,p'})^* \Bigr)  &\to 
 \Lambda^*_{\phi,p} \boxtimes (\Lambda^*_{\phi,p'})^*, \\
\bigl( (s, s''), (r'',r') \bigr) &\mapsto (s''(r'')s, s''(r'') r').
\end{split} \end{equation}

We will prove this theorem below, and assuming it for the moment 
we conclude our main result (Theorem \ref{theorem-main}) as a corollary.

\begin{Cor}\label{main-corr}
The heat kernel of $\Delta_\phi$ lifts to a polyhomogeneous section of $E$ in \eqref{E} on the
heat space $HM_\phi$, vanishing to infinite order at ff, tf, rf and lf, smooth at fd, and 
of order $(-m)$ at td. More precisely, $e^{-t\Delta_\phi} \in \mathscr{H}_\phi^{\, 3,0}(M, \Lambda^*_\phi)$.
Same holds for Laplace-type operators as in Remark \ref{dirac-squares}, 
with $\Lambda^*_\phi$ replaced by $L$. 
\end{Cor}

\begin{proof}
Consider the initial parametrix $H_0$ of Theorem \ref{initial-parametrix-th}.
When viewed as operators acting by an additional convolution in time, we find
$$
(\partial_t + \Delta_\phi) H_0 =   (\textup{Id} + P), \quad P \in \mathscr{H}_\phi^{\, 4,\infty}(M, \Lambda^*_\phi).
$$
Formally, the heat kernel is obtained by inverting the error term $(\textup{Id} + P)$
\begin{align}\label{Neumann0}
e^{-t\Delta_\phi} = H_0 \circ (\textup{Id} + P)^{-1} = \sum_{\ell = 0}^\infty (-1)^\ell H_0 \circ P^\ell.
\end{align}
However, convergence of the Neumann series in $\mathscr{H}_\phi^{\, *}(M, \Lambda^*_\phi)$ is intricate, since 
the lift of the Schwartz kernel for error term $P$ to $HM_\phi$ does not vanish to infinite order at fd.
We remedy this by correcting the fd asymptotics of $H_0$, asymptotically summing 
the lifts $\beta^*((-1)^{\ell} H_0 \circ P^{\ell})$ over $\ell \in \mathbb{N}_0$. This defines
a new heat parametrix $H'_0 \in \mathscr{H}_\phi^{\, 3,\infty}(M, \Lambda^*_\phi)$, such that 
\begin{equation}
\begin{split}
H'_0 \sim \sum_{\ell = 0}^\infty (-1)^\ell N_{\textup{fd}} \bigl( H_0 \circ P^\ell \bigr) \rho^{-3+4\ell}_{\textup{fd}}.
\end{split}
\end{equation}
As a consequence, we obtain a much better error term whose Schwartz kernel, lifted to $HM_\phi$, by construction 
vanishes to infinite order at all boundary faces of $HM_\phi$
\begin{equation}
\begin{split}
(\partial_t + \Delta_\phi) H'_0 =   (\textup{Id} + P'), \quad P' \in \mathscr{H}_\phi^{\, \infty, \infty}(M, \Lambda^*_\phi).
\end{split}
\end{equation}
Now we can invert $(\textup{Id} + P')$ and obtain the heat kernel as a Neumann series
\begin{align}\label{Neumann}
e^{-t\Delta_\phi} \equiv H = H'_0 \circ (\textup{Id} + P')^{-1} = H'_0 + \sum_{\ell = 1}^\infty (-1)^\ell H'_0 \circ P'^\ell.
\end{align}
The series can be shown to converge in $\mathscr{H}_\phi^{\, \infty, \infty}(M, \Lambda^*_\phi)$ 
by a Volterra series argument cf. \cite[Theorem 2.19]{berline},
and \cite[Proposition 2.10]{notes}. \medskip

We shall be precise: let $K_{P'}$ be the Schwartz kernel of $P'$,
and $K_{P'^\ell}$ the Schwartz kernel of $P'^\ell$. For any 
fixed $(p,p',t) \in \overline{M}^2 \times \R^+$ we set 
$(p_0,t_0) = (p,t)$ and $(p_\ell,t_\ell) =(p',0)$, and then write $K_{P'^\ell}(p,p',t)$ as follows
\begin{align}\label{composition-ell}
K_{P'^\ell}(p,p',t) = 
\int_{\triangle} \int_{M^{\ell - 1}} 
\prod\limits_{k=0}^{\ell-1}K_{P'}(p_k,p_{k+1},t_k-t_{k+1}) 
\prod\limits_{k=0}^{\ell-2} dt_{k+1} \textup{dvol}_{g_\phi} (p_{k+1}),
\end{align}
where the integration region $\triangle \in [0,t]^{\ell-1}$ is a simplex given by
$$
\triangle = \left\{ (t_1, \cdots, t_{\ell-1}) \in [0,t]^{\ell-1} \mid t \geq t_1 \geq \cdots \geq t_{\ell - 1}\right\}.
$$
The lift $\beta^*K_{P'}$  to $HM_\phi$ vanishes to 
infinite order at all boundary faces. Thus, for any $N \in \mathbb{N}$ there 
exists a constant $C_N>0$ depending only on $N$ such that (we denote by $\| \cdot \|_\infty$ the
supremum over $\partial M \times \partial M$ of the pointwise ($g_\phi$-induced) norms on the fibres of $\Lambda^*_\phi \boxtimes (\Lambda^*_\phi)^*$)
$$
\|  K_{P'}(x, x',t)\|_\infty \leq C_N (t \, x \, x')^N (x')^{2+b}.
$$
Recall now the notation of Definition \ref{def-phi-metric} and consider a Riemannian metric $\overline{g}$ on $M$, such that 
$\overline{g} \restriction \mathscr{U} = dx^2 + g_F + \phi^* g_B$. Note that $M$ has finite 
volume $\textup{vol} (M, \overline{g})$ with respect to $\overline{g}$.
Moreover, we obtain for the volume forms of $g_\phi$ and $\overline{g}$
(up to a bounded nowhere vanishing factor)
$$
\textup{dvol}_{g_\phi} = (x')^{-2-b} \textup{dvol}_{\overline{g}}.
$$
Consequently, we obtain from \eqref{composition-ell} for the pointwise ($g_\phi$-induced) norms on the fibres of 
$\Lambda^*_\phi \boxtimes (\Lambda^*_\phi)^*$, exactly as in \cite[Proposition 2.10]{notes}
\begin{equation}
\begin{split}
\|  K_{P'^\ell}(x,y,z,x',y',z',t)\| &\leq C_N (t \, x \, x')^N \textup{vol}(\triangle) \textup{vol} (M, \overline{g}))^{\ell-1}
\\ &= (t \, x \, x')^N \frac{(t \, C_N \textup{vol} (M, \overline{g}))^{\ell-1}}{(\ell-1)!}.
\end{split}
\end{equation}
Similarly, we obtain an estimate for the Schwartz kernel of $H'_0 \circ P'^\ell$
\begin{equation}
\begin{split}
\|  K_{H'_0 \circ P'^\ell}(x,y,z,x',y',z',t) \| \leq C (t \, x \, x')^K \frac{(t \, C_N \textup{vol} (M, \overline{g}))^{\ell-1}}{(\ell-1)!},
\end{split}
\end{equation}
for some constant $C>0$, depending only on $H'_0$. Consequently
\begin{equation}
\begin{split}
\|  \Bigl( K_{H} - K_{H'_0} \Bigr) (x,y,z,x',y',z',t) \| \leq C (t \, x \, x')^N e^{t \, C_N \textup{vol} (M, \overline{g})}.
\end{split}
\end{equation}
Since $K \in \mathbb{N}$ was arbitrary, this proves that the infinite sum in \eqref{Neumann} converges in $\mathscr{H}_\phi^{\, *}(M, \Lambda^*_\phi)$.
This concludes the proof.
\end{proof}

\subsection{Proof of the composition theorem} 

In this subsection we prove Theorem \ref{compth}, where $\beta^*K_B$ is assumed to be 
vanishing to infinite order at td, but $\beta^*K_A$ is not necessarily.
We write the composition integral \eqref{prop comp}
using pullback and push forward as follows. We write $\R_+\equiv \R^+ := [0,\infty)$ and define the maps 
\begin{align*}
\pi_{C}&: \overline{M}^3\times \mathbb{R}_{t'}^{+} \times
\mathbb{R}_{t''}^{+} \longrightarrow \overline{M}^2 \times
\mathbb{R}^+_{t'+t''},\quad
(p,p',p'',t',t'') \rightarrow
(p,p'',t'+t''),\\
\pi_{L}&: \overline{M}^3\times \mathbb{R}_{t'}^{+} \times
\mathbb{R}_{t''}^{+} \longrightarrow \overline{M}^2\times
\mathbb{R}^+_{t''},\quad
(p,p',p'',t',t'') \rightarrow
(p,p',t''),\\
\pi_{R}&: \overline{M}^3\times \mathbb{R}_{t''}^{+} \times
\mathbb{R}_{t''}^{+} \longrightarrow \overline{M}^2\times
\mathbb{R}^+_{t'},\quad
(p,p',p'',t',t'') \rightarrow
(p',p'',t').
\end{align*}
Then we can write \eqref{prop comp} 
by pulling back $K_{A}$ and $K_{B}$ to
$\overline{M}^3\times \mathbb{R}_{+}^{2}$ via $\pi_{L}, \pi_{R}$
and pushing forward (integrate) along $t= t'+t''$ and in $p'$ with
respect to $\textup{dvol}_{g_\phi}$
\begin{equation}\label{ABC}
K_{C} = (\pi_{C})_{*}\left( \frac{}{} \! (\pi_{L}^{*}K_{A})\cdot
(\pi_{R}^{*}K_{B})\right),
\end{equation}
where the dot is defined fibrewise as in \eqref{dot-def}. 
We prove the composition theorem by constructing 
the "heat triple space" $HM^3_{\phi}$ by a resolution process
from $\overline{M}^3\times \mathbb{R}_{+}^{2}$, with blow down map
\begin{equation*}
\beta_{3}: HM_{\phi}^{3} \longrightarrow \overline{M}^3 \times \mathbb{R}_{+}^{2}.
\end{equation*}
We show that the maps $\pi_{C},\pi_{L},\pi_{R}$ lift
to b-fibrations $\Pi_{C},\Pi_{L},\Pi_{R}$ in the sense of Melrose \cite{MelATP}, i.e.
in the commutative diagram diagram \eqref{eq:0.3}
\begin{equation}\label{eq:0.3}
\begin{tikzcd}
HM^3_{\phi}\arrow[r]\arrow[d, "\beta_{\textup{Tr}}" ']& M^2_{\phi} \times \R^+ \arrow[d, "\beta_\phi"] \\ 
\overline{M}^3\times\mathbb{R}^2_{+}\arrow[r]&
 \overline{M}^2\times \mathbb{R}_{+}
\end{tikzcd}
\end{equation}
the maps $\pi_{C},\pi_{L},\pi_{R}$ in bottom arrow  lift
to b-fibrations $\Pi_{C},\Pi_{L},\Pi_{R}$ in the upper arrow.
 Here we use the
intermediate heat space $M^2_{\phi}\times \R^+$, introduced in Figure \ref{intermediate-blowup},
together with the corresponding blowdown map $\beta_\phi$, since the additional blowup of the temporal diagonal 
will not be necessary due to Proposition \ref{prop comp}.
Defining $\kappa_{A,B,C} := \beta_\phi^{*}(K_{A,B,C})$,
we obtain using the commutativity of  diagram \eqref{eq:0.3}
a new version of the relation \eqref{ABC}
\begin{equation}\label{eq:0.2} 
\kappa_{C} \equiv \beta_\phi^*(K_{C}) = (\Pi_{C})_{*}(\Pi_{L}^{*}\kappa_{A}\cdot \Pi_{R}^{*}\kappa_{B}).
\end{equation}
The idea is now to see that $\Pi_{L}^{*}\kappa_{A} \cdot \Pi_{R}^{*}\kappa_{B}$ 
is indeed polyhomogeneous and the pushforward under $\Pi_{C}$
preserves the polyhomogeneity. In the rest of this section we follow this strategy more
concretely, first construct the triple space $HM_{\phi}^{3}$, 
compute the lift of boundary defining functions under projections to
compute explicitly the asymptotics of
$\Pi_{L}^{*}\kappa_{A}, \Pi_{R}^{*}\kappa_{B}$ and also the pushfoward,
 $(\Pi_{C})_*(\Pi_{L}^{*}\kappa_{A}\cdot \Pi_{R}^{*}\kappa_{B})$.

\subsubsection*{Construction of the triple space}

\begin{figure}[h]
\centering
\begin{tikzpicture}[scale =0.5]
\draw[-] (0,2)--(0,5);
\draw[-] (-2,-0.5)--(-4,-1.5);
\draw[-] (2,-0.5)--(4,-1.5);

\node at (0,0.4) {$111$};
\node at (-2.7,1) {$001$};           
\node at (2.7,1) {$100$};
\node at (0,-1.7) {$010$};
         
 \draw (0,2).. controls (-1.5,1.5) and (-1.6,0.5) .. (-2,-0.5);           
 \draw (0,2).. controls (1.5,1.5) and (1.6,0.5) .. (2,-0.5);
 \draw (-2,-0.5).. controls (-1.7,-1) and (1.7,-1) .. (2,-0.5);           
           
\begin{scope}[shift={(12,0)}, scale = 1]
       \draw[->](0,0)--(0,5);
       \draw[->](0,0)--(-3,-1.5);
       \draw[->](0,0)--(3,-1.5);
\node at (-3.2,-2.2) {$x$};
\node at (3.5,-2) {$x''$};
\node at (1,5) {$x^{'}$};
  
\end{scope}
\draw[->] (5,1.5)-- node[above] {$\beta_{1}$} (7,1.5);
\end{tikzpicture}
\caption{Illustration of $M^3_{b}$ in spatial direction.}
\label{triple first}
\end{figure}
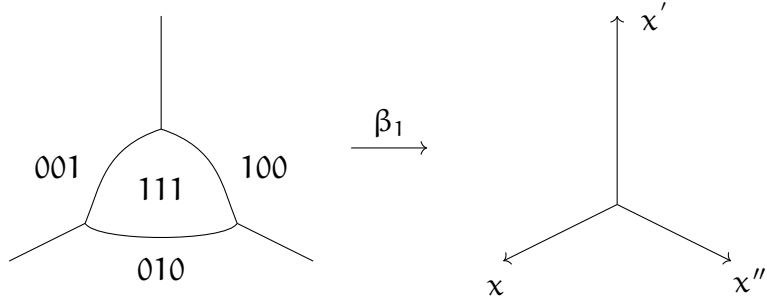

In order to apply Melrose's pushforward theorem \cite{Mel-push}, to 
conclude polyhomogeneity $\kappa_{C}$, the maps $\Pi_{C},\Pi_{L},\Pi_{R}$ need
to be b-fibrations. This dictates the construction of the 
triple space $HM_{\phi}^{3}$ as a blowup of $\overline{M}^3\times\mathbb{R}^2_{+}$. 
We describe the blowups using local coordinates $p= (x,y,z)$ and their 
copies $p'=(x',y',z')$ and $p''=(x'',y'',z'')$ on $M$. The time coordinates on each $\mathbb{R}^2_{+}$
are written as $t$ and $t'$. The first submanifold to blow up is then 

\begin{equation*}
F := \{t' = t'' = 0, x = x' = x'' = 0\} \subset \overline{M}^3\times\mathbb{R}^2_{+}. 
\end{equation*}
We refer the reader to \S \ref{elliptic-section} for the basic elements of the blowup procedure.
As before we blow up parabolically in the time direction, i.e. we treat 
$\sqrt{t}$ and $\sqrt{t'}$ as smooth. The resulting blowup space $M^3_{b} = [ \overline{M}^3\times \mathbb{R}^2_{+} ; F]$ 
is illustrated in Figure \ref{triple first} and comes with the blow down map

\begin{equation*}
\beta_{1} : M^3_{b} \longrightarrow \overline{M}^3\times \mathbb{R}^2_+.
\end{equation*}

We denote the resulting new boundary face, which is the inward spherical 
normal bundle of $F\subset \overline{M}^3\times\mathbb{R}^2_{+}$ by $(111)$,
where the first $1$ indicates that the boundary face corresponds to $x=0$, the second $1$ 
corresponds to $x'=0$, and the third $1$ to $x''=0$. This principle is also used in the namesgiving for other boundary
faces, e.g. $(100)$ is the lift of $\{x=0\}$. \medskip

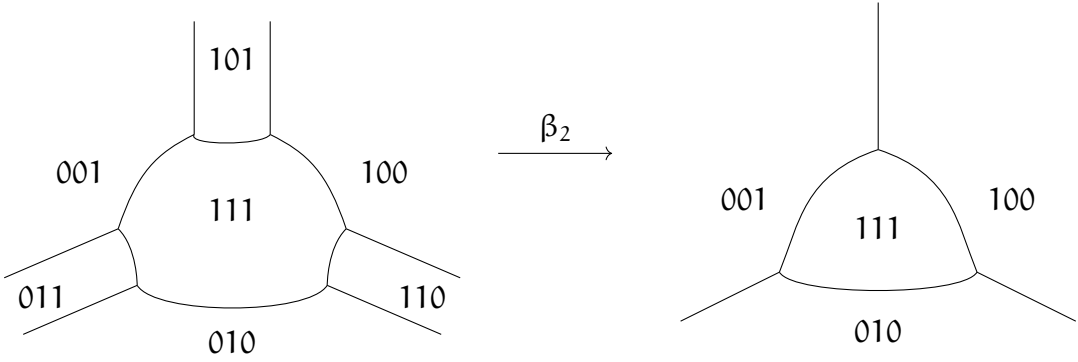
\begin{figure}[h]
\centering
\begin{tikzpicture}[scale =0.5]
      \draw[-](-1,2)--(-1,5);
       \draw[-](1,2)--(1,5);
       
       \draw[-](-3,-0.5)--(-6,-1.8);
       \draw[-](-2.5,-2)--(-5.5,-3.3);
       
           \draw[-](3,-0.5)--(6,-1.8);
       \draw[-](2.5,-2)--(5.5,-3.3);

\draw (-1,2).. controls (-1.2,1.7) and (0.8,1.7) .. (1,2);
\draw (-1,2).. controls (-2,1.5) and (-2.5,1) .. (-3,-0.5);
\draw (1,2).. controls (2,1.5) and (2.5,1) .. (3,-0.5);
    
\draw (-3,-0.5).. controls (-2.5,-1) and (-2.5,-2) .. (-2.5,-2);
\draw (3,-0.5).. controls (2.5,-1) and (2.5,-2) .. (2.5,-2);
\draw (-2.5,-2).. controls (-2,-2.8) and (2,-2.8) .. (2.5,-2);

 \node at (0,0) {$111$};
 \node at (0,-3.5) {$010$};
 \node at (0,4) {$101$};
 \node at (4,1) {$100$};
 \node at (-4,1) {$001$};
 \node at (5,-2.3) {$110$};
 \node at (-5,-2.3) {$011$};
                                       
\begin{scope}[shift={(17,-1)}, scale = 1.3]
\draw[-] (0,2)--(0,5);
\draw[-] (-2,-0.5)--(-4,-1.5);
\draw[-] (2,-0.5)--(4,-1.5);
\node at (0,0.4) {$111$};
           
  \node at (-2.7,1) {$001$};           
\node at (2.7,1) {$100$};
\node at (0,-1.7) {$010$};
       
 \draw (0,2).. controls (-1.5,1.5) and (-1.6,0.5) .. (-2,-0.5);           
 \draw (0,2).. controls (1.5,1.5) and (1.6,0.5) .. (2,-0.5);
 \draw (-2,-0.5).. controls (-1.7,-1) and (1.7,-1) .. (2,-0.5);

\end{scope}
\draw[->] (7,1.5)-- node[above] {$\beta_{2}$} (10,1.5);
\end{tikzpicture}
\caption{Illustration of $M^3_{b,t}$ in spatial direction.}
\label{triple second}
\end{figure}

\noindent We then blow up $M^3_{b}$ at 
$$ 
F_O := \beta^*_1 \{t'=t''=0\}.
$$
The next submanifolds to blow up are  submanifolds
of the codimension 2 in spatial direction
corresponding to each projection $\pi_{C}, \pi_{L}, \pi_{R}$.
Accordingly we denote these submanifolds as
 $F_{C}, F_{L}, F_{R}$, which are defined as
\begin{align*}
F_{C} &:= \beta^*_1 \{t' = t'' = 0, x = x'' = 0\},\\
F_{L} &:= \beta^*_1 \{t'' = 0, x' = x'' = 0\},\\
F_{R} &:= \beta^*_1 \{t' = 0, x = x' = 0\}.
\end{align*}
We point out that the order of blowing of submanifolds
$F_{C}, F_{L},F_{R}$ after blowing up of $F$ and $F_O$ is immaterial
as they become disjoint. As before we blow up parabolically in the time direction. 
The resulting blowup space 
$$
M^3_{b,t} =  [[M^3_{b}; F_O] ; F_{C},F_{L},F_{R}],
$$
is illustrated in Figure \ref{triple second} and comes with the blow down map
\begin{equation*}
\beta_{2} : M^3_{b,t} \longrightarrow M^3_{b}.
\end{equation*}

\begin{figure}[h]
  \includegraphics[width=\linewidth]{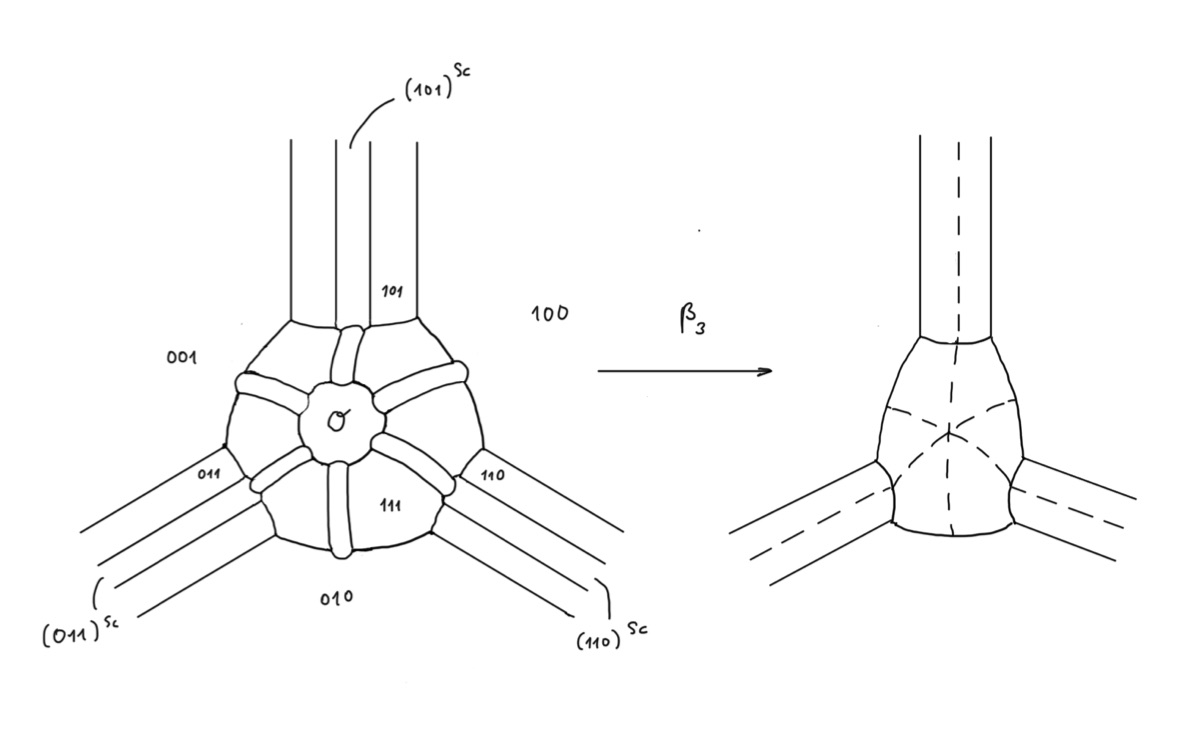}
  \caption{Illustration of $HM_{\phi}^{3}$ in spatial direction.}
  \label{triple third}
  \end{figure}

The triple elliptic space of the $\phi$-calculus, see Grieser and Hunsicker \cite{grieser-pseudo}, includes
the fibre-diagonal blow up in each direction. Here, we need to perform the same blowups
combined with blowing up the time direction. More precisely, using local coordinates, we blow up the following submanifolds
\begin{align*}
F_{\text{C,Sc}} &:= (F \cup F_C) \cap (\beta_2 \circ \beta_1)^* \{ x = x'' , y = y''\},\\
F_{\text{L,Sc}} &:= (F \cup F_L) \cap (\beta_2 \circ \beta_1)^* \{ x' = x'' , y' = y''\},\\
F_{\text{R,Sc}} &:= (F \cup F_R) \cap (\beta_2 \circ \beta_1)^* \{ x = x' , y = y'\}.
\end{align*}
as well as their intersection
\begin{equation*} 
O := F_{\text{C,Sc}} \cap 
F_{\text{L,Sc}} \cap
F_{\text{R,Sc}}.
\end{equation*}
This defines the triple space in the heat calculus, illustrated in Figure \ref{triple third} as
\begin{equation*}
HM_{\phi}^{3} := [\overline{M}^3_{b,t}, O; F_{\text{C,Sc}}, F_{\text{L,Sc}}, F_{\text{R,Sc}}].
\end{equation*}

The triple space comes with the intermediate blow down map 
$\beta_3: HM_{\phi}^{3} \to M^3_{b,t}$ and we define the total
blowdown map as
\begin{equation*}
\beta_{\textup{Tr}}:= \beta_3 \circ \beta_2 \circ \beta_1 : HM_{\phi}^{3} \longrightarrow \overline{M}^3\times \mathbb{R}^2_+.
\end{equation*}
From the construction, we compute immediately 
\begin{equation} \label{def triple}
\begin{split}
&\beta_{Tr}^{*}(x) = \rho_{111}\rho_{O}
\rho_{101}\rho_{101}^{sc}\rho_{110}\rho_{110}^{sc}\rho_{100},\\
&\beta_{Tr}^{*}(x') = \rho_{111}\rho_{O}
\rho_{011}\rho_{011}^{sc}\rho_{110}\rho_{110}^{sc}\rho_{010},\\
&\beta_{Tr}^{*}(x'') = \rho_{111}\rho_{O}
\rho_{101}\rho_{101}^{sc}\rho_{011}\rho_{011}^{sc}\rho_{001}.\\
\end{split}
\end{equation}

\subsubsection*{Lifts of boundary defining functions under various projections}
\begin{figure}[h]
  \includegraphics[scale=0.35]{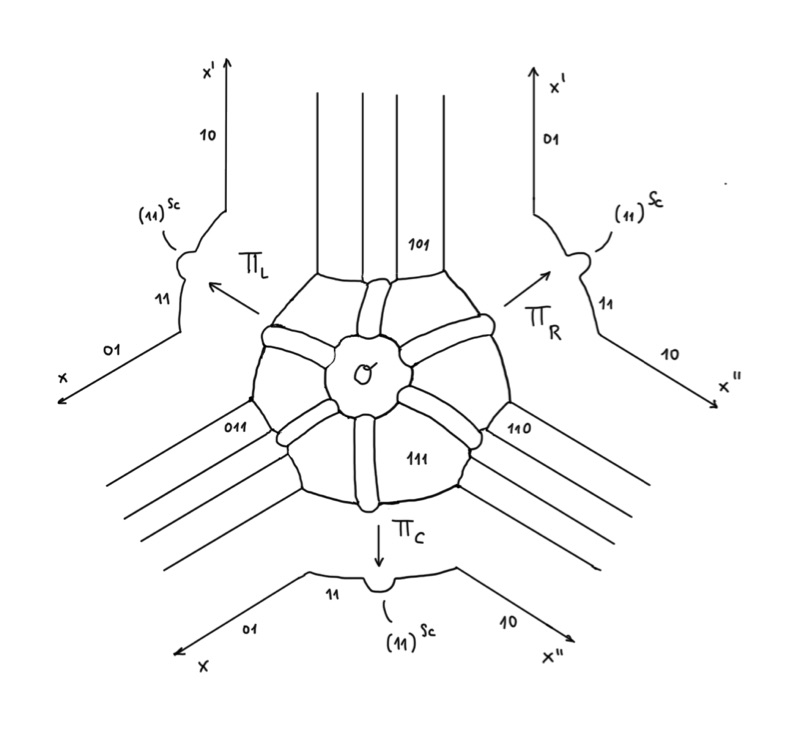}
  \caption{Illustration of projections in spatial direction.}
  \label{triple-proj}
  \end{figure}

We shall now study the lifts of boundary defining functions on the intermediate heat space $M^2_{\phi}\times \R^+$ 
to the triple heat space under the maps $\Pi_{C},\Pi_{L},\Pi_{R}: HM^3_{\phi} \to M^2_{\phi}\times \R^+$. We use the following 
notation: we denote a boundary defining function of any boundary face $(ijk)$, in $HM^3_{\phi}$ by 
$\rho_{ijk}; i,j,k \in \{0,1\}$. Boundary defining functions of $(110)^{\text{Sc}}, (101)^{\text{Sc}}$ and $(011)^{\text{Sc}}$ are 
denoted by $\rho_{110}^{\text{Sc}}, \rho_{101}^{\text{Sc}}$ and $\rho_{011}^{\text{Sc}}$, respectively. The 
boundary face $\mathcal{O}$, arising from the blowup of $O$, comes with a defining function $\rho_O$.
The boundary face, arising from the blowup of $F_O$, comes with a defining function $\tau_O$. 
Let $\tau$ and $\tau'$ be defining functions for the two boundary faces in $HM^3_{\phi}$, 
corresponding to $\{t'=0\}$ and $\{t''=0\}$, respectively. The notation is according to the labels in 
Figures \ref{triple-proj} and \ref{triple third}.
\medskip

We rename the boundary faces lf, rf and ff in the intermediate heat space $M^2_{\phi}\times \R^+$ as 
$(ij), i,j \in \{0,1\}$, according to the labels in Figure \ref{triple-proj} and write for the corresponding 
boundary defining functions $\rho_{ij}$. The boundary face fd is renamed $(11)^{\text{Sc}}$ and its 
defining function is written as $\rho_{11}^{\text{Sc}}$. Defining function of tf in  $M^2_{\phi}\times \R^+$ is denoted by $\tau$.
The notation is according to the labels in Figure \ref{triple-proj}, and corresponds closely to the notation of
boundary faces on the triple space. 

We compute  the pullback
of $\rho_{11},\rho_{01},\rho_{10},\rho_{11}^{\text{sc}}$
under  $\Pi_{C},\Pi_{L},\Pi_{R}$. Here, Figure \ref{triple-proj} provides a
helpful orientation.

\begin{align}\label{0.3}\begin{array}{lll}
\Pi_{\text{C}}^{*}(\rho_{11}) = \rho_{111}\cdot \rho_{101} &
 \Pi_{\text{L}}^{*}(\rho_{11}) = \rho_{111}\cdot \rho_{110} &
\Pi_{\text{R}}^{*}(\rho_{11}) = \rho_{111}\cdot \rho_{011}\\
\Pi_{\text{C}}^{*}(\rho_{01}) = \rho_{011}\cdot \rho_{011}^{\text{Sc}}\cdot \rho_{001} &
 \Pi_{\text{L}}^{*}(\rho_{01}) = \rho_{011}\cdot \rho_{011}^{\text{Sc}}\cdot \rho_{010} &
\Pi_{\text{R}}^{*}(\rho_{01}) = \rho_{101}\cdot \rho_{101}^{\text{Sc}}\cdot \rho_{001} \\
\Pi_{\text{C}}^{*}(\rho_{10}) = \rho_{110}\cdot \rho_{110}^{\text{Sc}}\cdot \rho_{100} &
 \Pi_{\text{L}}^{*}(\rho_{10}) = \rho_{101}\cdot \rho_{101}^{\text{Sc}}\cdot \rho_{100} &
 \Pi_{\text{R}}^{*}(\rho_{10}) = \rho_{110}\cdot \rho_{110}^{\text{Sc}}\cdot \rho_{010} \\
\Pi_{\text{C}}^{*}(\rho_{11}^{\text{Sc}}) = \rho_{101}^{\text{Sc}}\cdot \rho_O&
 \Pi_{\text{L}}^{*}(\rho_{11}^{\text{Sc}}) = \rho_{110}^{\text{Sc}}\cdot \rho_O&
 \Pi_{\text{R}}^{*}(\rho_{11}^{\text{Sc}}) = \rho_{011}^{\text{Sc}}\cdot \rho_O
\end{array}\end{align}

Now the lifts of the time variable $\tau$ are somewhat more intricate to argue.
Let us first compute the lifts of time direction boundary defining
functions under the blow down map $\beta_{\textup{Tr}}$. We find
\begin{equation} \label{eq:05}
\begin{split}
&\beta^*_{\textup{Tr}}(t') = \tau'\cdot\tau_O\cdot \rho_{110}^{2}
\cdot (\rho_{110}^{\text{Sc}})^{2}\cdot  \rho_{101}^{2}
\cdot (\rho_{101}^{\text{Sc}})^{2}\cdot \rho_{111}^{2}\cdot \rho_O^{2}, \\
&\beta^*_{\textup{Tr}}(t'') = \tau''\cdot\tau_O\cdot
 \rho_{101}^{2}\cdot (\rho_{101}^{\text{Sc}})^{2}\cdot
   \rho_{011}^{2}\cdot (\rho_{011}^{\text{Sc}})^{2}\cdot
    \rho_{111}^{2}\cdot \rho_O^{2},\\
&\beta^*_{\textup{Tr}}(t' + t'') = \tau_O\cdot
 \rho_{111}^{2}\cdot \rho_{101}^{2}\cdot
  \rho_{F_{0}}^{2} \cdot (\rho_{101}^{\text{Sc}})^{2}
\end{split}
\end{equation}
From the commutative diagram $\eqref{eq:0.3}$, we have
\begin{equation}\label{pb}
\pi_{C}\circ \beta_{\textup{Tr}} = \beta_\phi\circ\Pi_{C}.
\end{equation}
Now the lifts $\beta^*_\phi(t'), \beta^*_\phi(t''), \beta^*_\phi(t' + t'')$
to $\Pi_{L}(HM^3_{\phi}), \Pi_{R}(HM^3_{\phi}), \Pi_{C}(HM^3_{\phi})$,
respectively, are equal to
$\tau\cdot(\rho_{11}^{Sc})^{2}\cdot\rho_{11}^{2}$.
Therefore we compute, in view of $\eqref{eq:05}$ and \eqref{pb}
\begin{align*}
\Pi_{C}^{*}(\tau(\rho_{11}^{\text{Sc}})^{2}\rho_{11}^{2})
&= \Pi_{C}^{*}(\tau)\rho_{111}^{2}(\rho_{101}^{\text{Sc}})^{2}
\cdot\rho_O^{2},\\ &= \tau_{O}\cdot
\rho_{111}^{2}\cdot(\rho_{101}^{\text{Sc}})^{2}\cdot
\rho_O^{2}.
\end{align*}
From here we conclude $\Pi_{C}^{*}(\tau) = \tau_O$.
Similarly one can compute the other lifts and we arrive at the following 
identities
\begin{equation}\label{tlifts}
\begin{split}
\Pi_{L}^{*}(\tau) &= \tau'\cdot\tau_O\cdot(\rho_{101}^{\text{Sc}})^{2}\cdot\rho_{101}^{2},\\
\Pi_{R}^{*}(\tau) &= \tau'\cdot\tau_O\cdot(\rho_{101}^{\text{Sc}})^{2}\cdot\rho_{101}^{2}, \\
\Pi_{C}^{*}(\tau) &= \tau_O.
\end{split}\end{equation}

\subsubsection*{Projections $\Pi_{C},\Pi_{L},\Pi_{R}$ are b-fibrations}
Condition for the application of Melrose's pushforward theorem \cite{Mel-push}
is that the maps $\Pi_{C},\Pi_{L},\Pi_{R}$ are b-fibrations, recall Definition 
\ref{b-maps-def}. By discussion above, \eqref{0.3} and \eqref{tlifts},
$\Pi_{C},\Pi_{L},\Pi_{R}$ are indeed b-fibrations.

\subsubsection*{Lifts of kernels and densities to the triple space}
Consider $A \in \mathscr{H}_\phi^{\ell,q}(M, \Lambda^*_\phi), B \in \mathscr{H}_\phi^{\ell',\infty}(M, \Lambda^*_\phi)$. 
We write $K_A$ and $K_B$ for 
the Schwartz kernels of $A$ and $B$, respectively. We consider the composition 
$C= A \circ B$ with Schwartz kernel $K_C$. We have by construction
\begin{equation}
\begin{split}
&\pi^*_L K_A = K_{\text{A}}(t',x,y,z,x',y',z'), \\
&\pi^*_R K_B = K_{\text{B}}(t'',x',y',z',x'',y'',z''). 
\end{split}
\end{equation}
We also write $K_C \equiv K_{\text{C}}(t,x,y,z,x'',y'',z'')$ and set
\begin{equation}
\begin{split}
&\nu_3 := dt'dt''\textup{dvol}_{g_\phi}(x,y ,z) \, \textup{dvol}_{g_\phi}(x',y',z') \, \textup{dvol}_{g_\phi}(x'',y'',z''), \\
&\nu_2 := dt \, \textup{dvol}_{g_\phi}(x,y ,z) \, \textup{dvol}_{g_\phi}(x'',y'',z'').
\end{split}
\end{equation}
Then we obtain by construction (making the relation \eqref{ABC} precise)
\begin{equation*}
K_C \cdot  \nu_2 = (\pi_C)_* \left( \pi^*_L K_A \cdot \pi^*_R K_A \cdot \beta_{Tr}^{*} \nu_3 \right). 
\end{equation*}
Writing this relation in terms of the lifts $\kappa_{A,B,C} = \beta^*_\phi K_{A,B,C}$,
we obtain
\begin{equation}\label{ABC-lifted}
\kappa_C \cdot  \beta_\phi^* \nu_2 = (\Pi_C)_* \left( \Pi^*_L \kappa_A \cdot \Pi^*_R \kappa_A \cdot \beta^*_{\textup{Tr}} \nu_3 \right). 
\end{equation}
This formula makes clear how to proceed \medskip

\begin{enumerate}
\item[] Step 1):  Compute the asymptotics of $\Pi^*_L \kappa_A \cdot \Pi^*_R \kappa_A$.
\item[] Step 2): Compute the asymptotics of $\Pi_{L}^{*}\kappa_{A} \cdot \Pi_{R}^{*}\kappa_{B} \cdot \beta^*_{\textup{Tr}} \nu_3$.
\item[] Step 3): Apply Pushfoward Theorem to study pushforward by $(\Pi_C)_*$ 
\item[] Step 4): Compute the asymptotics of $\beta_\phi^* \nu_2$
\item[] Step 5): Compare the asymptotics of both sides in \eqref{ABC-lifted} to study $\kappa_C$.
\end{enumerate}

\subsubsection*{Step 1):  Compute the asymptotics of $\Pi^*_L \kappa_A \cdot \Pi^*_R \kappa_A$.}
Since $B \in \mathscr{H}_\phi^{\ell',\infty}(M, \Lambda^*_\phi)$, we find 
in view of the \eqref{0.3} and \eqref{tlifts} after an explicit
counting of exponents
\begin{align*}
\kappa_{B} \sim \rho_{\textup{fd}}^{ \ell'-3}
\left( \rho_{\textup{ff}} \rho_{\textup{rf}}  \rho_{\textup{lf}} \rho_{\textup{tf}} \right)^{\infty} 
\equiv  (\rho_{11}^{\textup{Sc}})^{\ell'-3}
\left( \tau \rho_{11} \rho_{10}  \rho_{01}\right)^{\infty}, \quad \Pi_{R}^{*}\kappa_{B} &\sim (\rho_{011}^{\text{Sc}} \rho_O)^{\ell'-3},
\end{align*}
where we adopted a convention that for kernels on the triple space $HM_{\phi}^{3}$,
we only write out asymptotics at those boundary faces, where the kernel is not 
vanishing to infinite order, suppressing the other boundary faces from the formula.
\medskip

We would like to write down a similar expansion for $A \in \mathscr{H}_\phi^{\ell,p}(M, \Lambda^*_\phi)$.
However, since in general $p \neq \infty$, the asymptotics of $\kappa_A$ is not uniform as $\tau \to 0$. Nevertheless, 
by \eqref{tlifts} the asymptotics of the lifts of $\tau$ under $\Pi_L$ and $\Pi_R$ coincide. Since
$\Pi_R^*\kappa_{B}$ contributes $\Pi_R(\tau^\infty)$, the composition $\Pi_{L}^{*}\kappa_{A} \cdot \Pi_{R}^{*}\kappa_{B}$
is still polyhomogeneous and vanishing to infinite order at the faces $\mathcal{O}, (101)^{\textup{Sc}}, (101)$
and the lift of $\{t''=0\}$. Thus, we may assume $p=\infty$ without loss of generality and obtain
\begin{align*}
\kappa_{A} \sim \rho_{\textup{fd}}^{ \ell - 3}
\left( \rho_{\textup{ff}} \rho_{\textup{rf}}  \rho_{\textup{lf}} \rho_{\textup{tf}} \right)^{\infty} 
\equiv  (\rho_{11}^{\textup{Sc}})^{\ell -3}
\left( \tau \rho_{11} \rho_{10}  \rho_{01}\right)^{\infty}, \quad 
\Pi_{L}^{*}\kappa_{A} &\sim (\rho_{110}^{\text{Sc}} \rho_O)^{\ell-3}.
\end{align*}
Consequently, we arrive under the notation which suppresses defining functions 
of infinite order
\begin{align}\label{ABO}
\Pi_{L}^{*}\kappa_{A} \cdot \Pi_{R}^{*}\kappa_{B} \sim \rho_O^{\ell + \ell'-6}.
\end{align}

\subsubsection*{Step 2):  Compute the asymptotics of $\Pi_{L}^{*}\kappa_{A} \cdot \Pi_{R}^{*}\kappa_{B} \cdot \beta^*_{\textup{Tr}} \nu_3$.} \ \medskip

\noindent The density $\nu_3$ is given in local coordinates by 
\begin{align*}
\nu_3 &= dt'dt'' \left( x^{-b-2} dx dy dz\right)
\left( x'^{-b-2} dx' dy' dz' \right) 
\left( x''^{-b-2} dx'' dy'' dz'' \right) \\
&= t' t'' (x x' x'')^{-1} \frac{dt'}{t'} \frac{dt''}{t''} \left( \frac{dx}{x} \frac{dy}{x^b} dz \right)
 \left( \frac{dx'}{x'} \frac{dy'}{x'^b} dz' \right)
 \left( \frac{dx''}{x''} \frac{dy''}{x''^b} dz'' \right) \\
 &=: t' t'' (x x' x'')^{-1} \nu'_3,
\end{align*}
up to a smooth bounded function. The lift $\beta^*_{\textup{Tr}} \nu'_3$ equals $h \cdot \nu_b^{(3)}$, where $\nu_b^{(3)}$
is a b-density on $HM^3_{\phi}$,
that is a smooth density on $HM^3_{\phi}$, divided by a product of all its boundary
defining functions. The factor $h$ is a polyhomogeneous function on $HM^3_{\phi}$, 
smooth at the boundary face $\mathcal{O}$; its asymptotics at other boundary faces is irrelevant, 
since $\Pi_{L}^{*}\kappa_{A} \cdot \Pi_{R}^{*}\kappa_{B}$ vanishes to infinite order there. 
Thus we compute in view of \eqref{def triple}, \eqref{eq:05} and \eqref{ABO}
\begin{equation}\label{ll5}
\begin{split}
&\Pi_{L}^{*}\kappa_{A} \cdot \Pi_{R}^{*}\kappa_{B} \cdot \beta^*_{\textup{Tr}} \nu_3 
\\ &= \Pi_{L}^{*}\kappa_{A} \cdot \Pi_{R}^{*}\kappa_{B} \cdot \beta^*_{\textup{Tr}} \left( t' t'' (x x' x'')^{-1} \right) \beta^*_{\textup{Tr}} \nu'_3 \\
&\sim \rho_O^{\ell + \ell'-5} \beta^*_{\textup{Tr}} \nu'_3 \sim  \rho_O^{\ell + \ell'-5} \nu_b^{(3)},
\end{split}
\end{equation}
where as before we suppressed the other boundary defining functions from the notation, where
the kernels vanish to infinite order.

\subsubsection*{Step 3):  Apply Pushfoward Theorem to
 study pushforward by $(\Pi_C)_*$.} \ \medskip

\noindent The next step is applying the pushforward theorem of Melrose
\cite{Mel-push}. We use the notation of Definition \ref{b-maps-def}.
Then the pushforward theorem of Melrose says the following. 

\begin{Th}
Let $M,M'$ be two compact manifolds with corners and $\nu_b, \nu'_b$
are b-densities on $M,M'$, respectively. Let $u$ be a polyhomogeneous function on $M$ with
index sets $E_{j}$ at the faces $H_{j}$ of $M$.
Suppose that each $(z,p) \in E_{j}$ has $\text{Re} (z) > 0$
if the index $j$ satisfies $e(i,j) = 0$ for all $j$.
Then the pushforward $f_{*}(u\nu_{b})$ is well-defined
and equals $h\nu_{b}'$ where $h$ is polyhomogeneous on
$M'$ and has an index family $f_{b}(\mathcal{E})$
given by an explicit formula in terms of the
index family $\mathcal{E}$ for $M$.
\end{Th}

We refer the reader to \cite{Mel-push, MelATP} for the explicit definition of the index family $f_{b}(\mathcal{E})$,
and just say that in our specific case we have for any $\alpha > 0$
(all other suppressed boundary functions enter with infinite order)
$$
(\Pi_C)_* \left( \rho_O^{\alpha} \nu_b^{(3)}\right) = (\rho_{11}^{\text{Sc}})^{\alpha} \nu_b^{(2)},
$$
where $\nu_b^{(2)}$ is a b-density on the intermediate heat space $M^2_{\phi}\times \R^+$, that is
a smooth density on $M^2_{\phi}\times \R^+$, divided by a product of all its boundary
defining functions. Hence we arrive in view of \eqref{ll5} at
\begin{equation}\label{push1}
\begin{split}
(\Pi_C)_* \left(  \Pi_{L}^{*}\kappa_{A} \cdot \Pi_{R}^{*}\kappa_{B} \cdot \beta^*_{\textup{Tr}} \nu_3  \right) &= 
(\rho_{11}^{\text{Sc}})^{\ell + \ell'-5} \nu_b^{(2)} \\
&\equiv (\rho_{11}^{\text{Sc}})^{\ell + \ell'-5} \left( \tau \rho_{11} \rho_{10}  \rho_{01}\right)^{\infty}  \nu_b^{(2)}.
\end{split}
\end{equation}

\subsubsection*{Step 4):  Compute the asymptotics of $\beta_\phi^* \nu_2$.} \ \medskip

\noindent The density $\nu_2$ is given in local coordinates by 
\begin{align*}
\nu_2 &= dt \left( x^{-b-2} dx dy dz\right)
\left( x''^{-b-2} dx'' dy'' dz'' \right) \\
&= t' t'' (x x'')^{-1}  \frac{dt}{t} \left( \frac{dx}{x} \frac{dy}{x^b} dz \right)
 \left( \frac{dx''}{x''} \frac{dy''}{x''^b} dz'' \right) \\
 &=: t (x x'')^{-1} \nu'_2,
\end{align*}
up to a smooth bounded function. The lift $\beta^*_{\phi} \nu'_2$ equals $h' \cdot \nu_b^{(2)}$, where 
$h'$ is a polyhomogeneous function on $M^2_{\phi}\times \R^+$, 
smooth at the boundary face fd; its asymptotics at other boundary faces is irrelevant, 
since it will be multiplied with a kernel that vanishes to infinite order there. 
Hence we find 
\begin{equation}\label{push2}
\begin{split}
\beta_\phi^* \nu_2 &= \beta_\phi^* \left(t (x x'')^{-1} \right) \nu_b^{(2)} \\
&=  \tau \left(\rho_{11}^{\text{Sc}} \rho_{11}\right)^{-2} \left(\rho_{10}  \rho_{01}\right)^{-1} h' \nu_b^{(2)}
\end{split}
\end{equation}

\subsubsection*{Step 5):  Compare the asymptotics of both sides in \eqref{ABC-lifted} to study $\kappa_C$.} \ \medskip

\noindent Combining \eqref{push1} and \eqref{push2}, we find

\begin{equation}\label{push3}
\begin{split}
(\Pi_C)_* \left(  \Pi_{L}^{*}\kappa_{A} \cdot \Pi_{R}^{*}\kappa_{B} \cdot \beta^*_{\textup{Tr}} \nu_3  \right) &= 
(\rho_{11}^{\text{Sc}})^{\ell + \ell'-5} \left( \tau \rho_{11} \rho_{10}  \rho_{01}\right)^{\infty}  \nu_b^{(2)} \\ &=
(\rho_{11}^{\text{Sc}})^{\ell + \ell'-3} \left( \tau \rho_{11} \rho_{10}  \rho_{01}\right)^{\infty}  \beta_\phi^* \nu_2 \\ &= 
\kappa_C \cdot  \beta_\phi^* \nu_2.
\end{split}
\end{equation}
This proves the composition Theorem \ref{compth}.

\end{document}